\documentclass[aos,nameyear,dvips]{arximspdf}
\usepackage{graphics}
\doi{10.1214/10-AOS797}
\volume{38}
\issue{5}
\pubyear{2010}
\firstpage{2678}
\lastpage{2722}

\makeatletter

\newtheorem{lemma}{Lemma}
\newtheorem{corollary}[lemma]{Corollary}
\newproclaim{example}[lemma]{Example}
\newtheorem{proposition}[lemma]{Proposition}
\newtheorem{theorem}[lemma]{Theorem}

\makeatother

\begin{document}
\begin{frontmatter}

\title{Generalized density clustering\protect\thanksref{T1}}
\runtitle{Density clustering}

\thankstext{T1}{Supported by NSF Grants DMS-07-07059, DMS-06-31589,
ONR Grant N00014-08-1-0673 and
a Health Research Formula Fund Award granted by the Commonwealth of
Pennsylvanias Department of
Health.}

\begin{aug}
\author[A]{\fnms{Alessandro} \snm{Rinaldo}\corref{}\ead[label=e1]{arinaldo@stat.cmu.edu}} and
\author[A]{\fnms{Larry} \snm{Wasserman}\ead[label=e2]{larry@stat.cmu.edu}}
\runauthor{A. Rinaldo and L. Wasserman}
\affiliation{Carnegie Mellon University}
\address[A]{Department of Statistics\\
Carnegie Mellon University\\
Pittsburgh, Pennsylvania, 15213\\
USA\\
\printead{e1}\\
\phantom{E-mail: }\printead*{e2}} 
\end{aug}

\received{\smonth{7} \syear{2009}}
\revised{\smonth{11} \syear{2009}}

%
\begin{abstract}
We study generalized density-based clustering in which sharply defined
clusters such as clusters on lower-dimensional manifolds are allowed.
We show that accurate clustering is possible even in high dimensions.
We propose two data-based methods for choosing the bandwidth and we
study the stability properties of density clusters. We show that a
simple graph-based algorithm successfully approximates the high density
clusters.
\end{abstract}

%
\begin{keyword}[class=AMS]
\kwd[Primary ]{62H30}
\kwd[; secondary ]{62G07}.
\end{keyword}
\begin{keyword}
\kwd{Density clustering}
\kwd{kernel density estimation}.
\end{keyword}

\end{frontmatter}

\section{Introduction}

It has been observed that classification methods
can be very accurate in high-dimensional problems,
apparently contradicting the curse of dimensionality.
A plausible explanation for this phenomenon
is the ``low noise'' condition described, for instance, in \citet
{mammen1999}.
When the low noise condition holds,
the probability mass near the decision boundary is small
and fast rates of convergence of the
classification error are possible in high
dimensions.

Similarly,
clustering methods can be
very accurate in high-dimensional problems.
For example, clustering subjects based on gene profiles
and clustering curves are both high-dimensional problems
where several methods have worked well despite the
high dimensionality.
This suggests that it may be possible to find
conditions that explains the success of clustering
in high-dimensional problems.

In this paper, we focus on
clusters that are defined
as the connected components of high density regions
[\citet{cuevasfraiman1997}, \citet{hartigan}].
The advantage of density clustering over other methods
is that there is a well-defined population
quantity being estimated
and density clustering allows the shape of the clusters to be
very general.
(A related but somewhat different approach for generally shaped
clusters is
spectral clustering; see
\citet{Lux} and [\citet{jordan}].)
Of course, without some conditions,
density estimation is subject to the usual
curse of dimensionality.
One would hope that an appropriate low noise
condition would obviate
the curse of dimensionality.
Such assumptions have been proposed by
\citet{polonik1995},
\citet{Rigollet06},
\citet{rigvert2006},
and others.
However, the assumptions used by these authors
rule out the case where the clusters are very sharply defined, which
should be the easiest cases, and, more generally,
clusters defined on lower dimensional sets.

The purpose of this paper is
to define a notion of density clusters
that does not rule out the most favorable cases and is not limited to
sets of full dimension.
We study the risk properties of density-based clustering and its
stability properties,
and we provide data-based methods for choosing the
smoothing parameters.

The following simple example helps to illustrate our motivation. We
refer the reader to the next section for a more rigorous introduction.
Suppose that a distribution $P$ is a mixture of finitely many point
masses at distinct points $x_1,\ldots, x_k$ where $x_j\in
\mathbb{R}^d$. Specifically,\vspace*{1pt} suppose that $P = k^{-1}\sum_{j=1}^k
\delta_j$ where $\delta_j$ is a point mass at $x_j$. The clusters are
$C_1=\{x_1\}, \ldots, C_k=\{x_k\}$. This is a trivial clustering
problem even if the dimension $d$ is very high. The clusters could not
be more sharply defined yet the density does not even exist in the
usual sense. This makes it clear that common assumptions about the
density such as smoothness or even boundedness are not well-suited for
density clustering.

Now let $p_h=dP_h/d\mu$ be the Lebesgue density of the measure $P_h$
obtained by convolving $P$ with the probability measure having
Lebesgue density $K_h$, a kernel with bandwidth $h$. Unlike the
original distribution $P$, $P_h$
has full-dimensional support for each positive $h$. The ``mollified''
density $p_h$ contains
all the information needed for clustering. Indeed, there exist
constants $\overline{h} > 0$ and $\lambda\geq0$ such that the
following facts are true:
\begin{enumerate}
\item
for all
$0 < h < \overline{h}$,
the level set
$\{x\dvtx p_h(x) \geq\lambda\}$
has disjoint, connected components
$C_1^h, \ldots, C_k^h$;
\item
the components $C_j^h$ contain the true clusters:
$C_j \subset C_j^h$ for $j=1,\ldots, k$;
\item although $C_j^h$ overestimates the true cluster $C_j$,
this overestimation is inconsequential since
$P(C_j^h -C_j)=0$ and hence a new observation will not be misclustered;
\item
let $\widehat{p}_h$ denote the kernel density estimator using $K_h$
with fixed bandwidth
$0 < h < \overline{h}$ and based on a i.i.d. sample of size $n$ from $P$.
Then,
${\sup_x}|p_h(x) - \widehat{p}_h(x)| = O(\sqrt{\log n/n})$ almost
everywhere $P$,
which does not depend on the dimension $d$ (see Section \ref
{sec::preliminaries}).
The bias from using a fixed bandwidth $h$---which does not vanish as
$n \to\infty$---does not adversely affect the clustering.
\end{enumerate}

In summary,
we can recover the true clusters
using an estimator of the density $p_h$
with a large bandwidth $h$.
It is not necessary to assume that the true
density is smooth or that it even exists.

Our contributions in this paper are the following:
\begin{enumerate}
\item We develop a notion
of density clustering that
applies to probability distributions that have nonsmooth Lebesgue
densities or do not even admit a density.
\item We find the rates of convergence for estimators of these clusters.
\item We study two data-driven methods for choosing the bandwidth.
\item We study the stability properties of density clusters.
\item We show that the depth-first search algorithm
on the $\rho$-nearest
neighborhood graph of $\{ \widehat{p}_h \geq\lambda\}$ is effective
at recovering the high-density clusters.
\end{enumerate}

Another approach to clustering that does not require densities
is the minimum volume set approach
[\citet{polonik1995}, \citet{scott}].
Our approach is different because we are specifically
trying to capture the idea that kernel density estimates are useful
for clustering even when the density may not exist.

Section \ref{section::clusters}
contains notation and definitions.
Section \ref{section::rates} contains results on
rates of convergence.
We give a data-driven method for choosing the bandwidth in
Section~\ref{sec::bandwidth}.
Section \ref{section::stability}
contains results on cluster stability.
The validity of the graph-based algorithm
for approximating the clusters is proved in Section
\ref{sec::friend}.
Section \ref{section::examples} contains some examples based on
simulated data.
Concluding remarks are in Section \ref{section::discussion}.
All proofs are in the Section \ref{section::proofs}.
Some technical details are in the \hyperref[app]{Appendix}.

\subsection*{Notation}
For two sequences $ \{ a_n \}$ and $\{ b_n \}$, we write $a_n = O(b_n)$
and $a_n = \Omega(b_n)$ if there exists a constant $C>0$ such that,
for all $n$ large enough, $|a_n|/b_n \leq C$ and $|a_n|/b_n \geq C$,
respectively. If $a_n = \Omega(b_n)$ and $a_n = O(b_n)$, then we will
write $a_n \asymp b_n$. We denote with $\mathbb{P}(E)$ the probability
of a generic event $E$, whenever the underlying probability measure is
implicitly understood from the context. 
By the dimension of a Euclidean set, we will always mean the
$k$-dimensional Hausdorff dimension for some
integer $0\leq k \leq d$ (see the \hyperref[app]{Appendix}). These
sets may consist, for example, of smooth submanifolds or even single
points.

\section{Settings and assumptions}
\label{section::clusters}


\subsection{Level set clusters}\label{sec:settings}
In this section, we develop a probabilistic framework for the
definition of clusters we have adopted. For ease of readability,
the more technical measure-theoretic details are given in the \hyperref
[app]{Appendix}.

Let $P$ be a probability distribution on $\mathbb{R}^d$ whose support
$S$ (the smallest closed set of $P$-measure $1$) is comprised of an
unknown number $m$ of disjoint compact sets
$\{S_1,\ldots, S_m \}$ of different dimensions.
We define the \textit{geometric density} of $P$ as the
measurable function $p \dvtx\mathbb{R}^d \mapsto\mathbb{R}$ given
by
%
%
\begin{equation}\label{eq:geom.dens}
p(x) = \lim_{h \downarrow0}\frac{ P(B(x,h))}{v_d h^d},
\end{equation}
where $B(x,\varepsilon)$ is the Euclidean ball of radius $h$ centered at
$x$, $\mu$
is the $d$-dimen\-sional Lebesgue measure and $v_d \equiv\mu(B(0,1))$. Note
that, almost everywhere $P$, $p(x) = \infty$ if and only if $x$
belongs to some set $S_i$
having dimension strictly less than $d$ and is positive and finite if
and only if $x$ belongs to some $d$-dimensional set~$S_i$. In
general, $\int_{\mathbb{R}^d} p(x) \,d\mu(x) \leq1$ and, therefore, $p$
is not necessarily a probability density. Nonetheless, $p$ can be used
to recover
the support of $P$, since
\[
S = \overline{\{ x \dvtx p(x) > 0\}},
\]
where for a set $A \subset\mathbb{R}^d$, $\overline{A}$ denotes its
closure.

For $\lambda\geq0$, define the $\lambda$-level set
%
%
\begin{equation}
L \equiv L(\lambda) = \overline{ \{ x\dvtx p(x) \geq\lambda
\}},
\end{equation}
and its boundary $\partial L(\lambda) = \{ x \dvtx p(x) = \lambda\}$.
Throughout the paper, we will suppose that we are given a fixed value of
$\lambda< \| p \|_{\infty}$, where $\| p \|_{\infty} \equiv\sup_{x
\in\mathbb{R}^d}p(x)$.
Often, $\lambda$ is chosen so that
$P(L(\lambda)) \approx1-\alpha$ for some given $\alpha$.
In practice, it is advisable to present the results
for a variety of values of $\lambda$
as we discuss in Section \ref{section::discussion}.

We assume that there are $k \geq1$ disjoint, compact, connected sets
$C_1(\lambda),\ldots,\break C_k(\lambda)$ such that
\[
L = C_1(\lambda) \cup\cdots\cup C_k(\lambda).
\]
We will often write $C_j$ instead of of $C_j(\lambda)$ when the
dependence of $\lambda$ is clear from the context.
The value of $k$ is not assumed to be known. The sets $C_1, \ldots,
C_k$ are
called the \textit{$\lambda$-clusters} of $p$, or just
\textit{clusters}. In our setting, the $C_j$'s need not be full
dimensional. Indeed, $C_j$ might be a lower-dimensional manifold or
even a single point. Furthermore, if $S_i$ has dimension
smaller than $d$, then $C_j = S_i$, for some $j =1,\ldots,k$. Thus,
for any $\lambda\geq0$, the $\lambda$-clusters of $p$ will include
all the lower-dimensional components of $S$. On the other hand, if
$S_i$ is full-dimensional, then there may be multiple clusters
in it, depending on the value of $\lambda$.

We observe an i.i.d. sample $X = (X_1,\ldots, X_n)$
from $P$, from which we construct the kernel density estimator
%
%
\begin{equation}
\widehat{p}_h(x) = \frac{1}{n}\sum_{i=1}^n
\frac{1}{c_d h^d} K \biggl(\frac{x-X_i}{h} \biggr)\qquad \forall x
\in\mathbb{R}^d,
\end{equation}
where $c_d \equiv\int_{\mathbb{R}^d} K(x) \,d \mu(x)$. For simplicity,
we assume that the kernel $K \dvtx\mathbb{R}^d \mapsto\mathbb{R}_+$
is a symmetric, bounded, smooth function supported on the Euclidean
unit ball. These assumptions can be easily relaxed to include, for
instance, the case of regular kernels as defined in
Devroye, Gy\"{o}rfi and Lugosi [(\citeyear{dgl97}), Chapter 10].
In particular, while the compactness of the support of $K$
simplifies our analysis, it is not essential and could be replace by
assuming fast decaying tails. Further conditions on the kernel $K$ are
discussed in Section \ref{sec::preliminaries}.

Let $p_h \dvtx\mathbb{R}^d \mapsto\mathbb{R}$ be the measurable
function given by
%
%
\begin{equation}
p_h(x) = \int_S K_h(x - y) \,dP(y) = \mathbb{E} (\widehat{p}_h(x)) ,
\end{equation}
where $K_h(x) \equiv\frac{1}{c_d h^d} K (\frac{\Vert x\Vert}{h}
)$. Also, let $K_h \mu$ be
the probability measure given by $K_h \mu(A) = \int_A K_h(x) \,d
\mu(x)$, for any Borel set $A \subseteq\mathbb{R}^d$. Then, $p_h$
is the Lebesgue density of the
probability measure $P_h$ obtained by convolving $P$ with $K_h
\mu$. More precisely, for each measurable set $A$,
\[
P_h(A) = \int_{A} \int_S K_h(x - y) \,dP(y) \,d \mu(x) = \int_A p_h(x)
\,d \mu(x).
\]
Borrowing some terminology from analysis, where the
kernel $K$ is referred to as a mollifier, we call the measure $P_h$
and the density $p_h$ as the mollified measure and mollified density,
respectively. For each $h$, the mollification of $P$ by $K$ yields
that:
\begin{enumerate}
\item the mollified measure $P_h$ has full-dimensional support $S
\oplus B(0,h)$
and is absolutely continuous with respect to $\mu$;
here, for two set $A$ and $B$ in
$\mathbb{R}^d$, $A \oplus B \equiv\{ x + y \dvtx x \in A, y \in B\}
$ denotes its Minkowski sum;
\item the mollified density $p_h$ is of class $\mathcal{C}^\alpha$
whenever $K$ is of class $\mathcal{C}^\alpha$,
with $\alpha\in\mathbb{N}_+ \cup\{ \infty\}$.
(A real valued function is of class $\mathcal{C}^\alpha$ if its
partial derivatives up to order $\alpha$ exist and are continuous.)
\end{enumerate}
(As a referee pointed out to us, the properties of mollified measures
are related to the classical theory of distributions.)
Mollifying $P$ makes it better behaved. At the same time, $P_h$
and $p_h$ can be seen as approximations of the original measure $P$ and
the geometric density $p$, respectively, in a sense made precise by the
following result.
\begin{lemma}\label{lem:ph}
As $h \to0$, $P_h$ converges weakly to $P$ and
$
\lim_{h \to0} p_h(x) = p(x),
$
almost everywhere $P$.
\end{lemma}

To estimate the $\lambda$-clusters of $p$, we use the connected
components of $\widehat{L}$, that is, the $\lambda$-clusters of
$\widehat{p}_h$.
That is, we estimate $L$ with
%
%
\begin{equation}
\widehat{L}\equiv\widehat{L}_h(\lambda) = \{x \dvtx \widehat
{p}_h(x) \geq\lambda\}.
\end{equation}

In practice, finding the estimated clusters
is computationally difficult.
Indeed, to verify that two points
$x_1$ and $x_2$ are in the same cluster, we need to find at least one
path $\gamma\subset\mathbb{R}^d$ connecting them such that $\widehat
{p}_h(x) \geq\lambda$ for each $x \in\gamma$. Conversely, when
$x_1$ and $x_2$ do not belong to the same cluster, this property has to
be shown to fail for every possible path between them.
We discuss an algorithm for approximating
the clusters in Section \ref{sec::friend}.
Until then, we ignore the computational problems
and assume that the $\lambda$-clusters of $\widehat{p}_h$ can be
computed exactly.

\subsection{Risk}\label{sec:risk}

We consider two different risk functions.

\begin{itemize}
\item The \textit{level set risk} is defined to be
$R^L(p,\widehat{p}_h) = \mathbb{E}(\rho(p,\widehat{p}_h,P))$, where
%
%
\begin{equation}
\rho(r,q,P)=
\int_{ \{x \dvtx r(x) \geq\lambda\}\Delta\{x \dvtx q(x) \geq
\lambda\}} dP(x),
\end{equation}
and
$A \Delta B = (A \cap B^c)\cup(A^c \cap B)$
is the symmetric set difference.
\item Define the \textit{excess mass functional} as
%
%
\begin{equation}\label{eq:excess.mass}
\mathcal{E}(A) = P(A) - \lambda\mu(A)
\end{equation}
for any measurable set $A \subset\mathbb{R}^d$.
This functional is maximized by
the true level set $L$;
see
\citet{mueller1991} and \citet{polonik1995}.
We can use the excess mass functional as a risk function
except, of course, that we maximize it rather than minimize it.
Given an estimate $\widehat{L}$ of $L$ based on $\widehat{p}_h$,
we will then be interested in making the \textit{excess mass risk}
%
%
\begin{equation}
R^M(p,\widehat{p}_h) = \mathcal{E}(L) - \mathbb{E} (\mathcal
{E}(\widehat{L}) )
\end{equation}
as small as possible.
Furthermore, if $P$ has full-dimensional support, simple algebra
reveals that maximizing
$\mathcal{E}(A)$ is equivalent to minimizing,
\[
{\int_{A\Delta L}} |p - \lambda| \,d \mu,
\]
which is the loss function used by
\citet{nowak2007}. In this case, the minimizer $L$ is unique.
More generally, if
$P=P_0+P_1$ where
$P_0$ is the part of $P$ that is absolutely continuous with respect to the
Lebesgue measure, then
%
%
\begin{equation}\label{eq:excessmass}
\mathcal{E}(L) - \mathcal{E}(A) =
{\int_{A\Delta L} }|p_0 - \lambda| \,d \mu+ P_1(L) - P_1(A),
\end{equation}
where $p_0 = \frac{d P_0}{d \mu}$. It is clear that $L$ is no longer
the unique minimizer of the excess mass functional.

\end{itemize}

\subsection{Assumptions}

Throughout our analysis, we assume the following conditions.

\begin{enumerate}[(C2)]
\item[(C1)]
There exist positive constants $\gamma$, $C_1$ and $\overline
{\varepsilon}$ such that
\[
\mathbb{P} \bigl( | p(X) - \lambda| < \varepsilon\bigr)
\leq C_1 \varepsilon^\gamma\qquad \forall\varepsilon\in[0,\overline
{\varepsilon}).
\]
\item[(C2)]
There exist a positive constant $\overline{h}$,
and a permutation $\sigma$ of $\{1,\ldots,k\}$
such that, for all $h \in(0, \overline{h})$ and all $\lambda' \in
(\lambda- \overline{\varepsilon}, \lambda+ \overline{\varepsilon})$,
\[
L_{h}(\lambda') = \bigcup_{j=1}^k C_j^h(\lambda'),
\]
where:
\begin{enumerate}[(b)]
\item[(a)] $C_i^h(\lambda') \cap C_j^h(\lambda)' = \varnothing$ for
$1\leq i< j \leq k$;
\item[(b)] $C_j(\lambda') \subseteq C^h_{\sigma(j)}(\lambda')$, for
all $1 \leq j \leq k$.
\end{enumerate}
\item[(C3)] There exist a positive constant $C_2$ such that,
for all $h \in(0,\overline{h})$ and $\lambda' \in(\lambda-
\overline{\varepsilon}, \lambda]$,
\[
L(\lambda') = \bigcup_{j=1}^k C_j(\lambda'),
\]
where
\[
\mu\bigl( \partial C_j(\lambda') \oplus B(0,h)\bigr) \leq C_2 h^{(d -
d_i)\vee1},
\]
and $d_i$ is the dimension of the component $S_i$ of the support of $P$
such that $C_j(\lambda') \subseteq S_i$.
%
\end{enumerate}

\subsection{Remarks on the assumptions}

Conditions of the form (C1) or of other equivalent forms, are also
known as low noise condition or margin conditions in the classification
literature. They have appeared in many places, such as \citet
{tsybakov97}, \citet{mammen1999}, \citet
{baillocuestacuevas01}, \citet{tsybakov04}, \citet
{Steinwart05}, \citet{cuevasgonzalesrodrigez06}, \citet
{cadre06}, \citet{audiberttsybakov2007}, \citet
{castronowak08} and \citet{aarti09}.

This condition, first introduced in \citet{polonik1995},
provides a way to relate the stochastic fluctuations of
$\widehat{p}_h$ around its mean $p_h$ to the clustering risk. Indeed,
the larger $\gamma$, the smaller the effects of these
fluctuations, and the easier it is to obtain good clusters from
noisy estimates of $p_h$, for any $h < \overline{h}$.

Conditions (C2) simply require that the level set of the mollified
density include the true clusters. The additional fringe $L_h - L$ can
be viewed as a form of clustering bias. Though mild and reasonable,
these assumptions are particularly
important, as they imply that the estimated density
$\widehat{p}_h$ can be used quite effectively for clustering purposes,
for a range of bandwidth values. This is is shown in the next
simple result. Let $N(\lambda)$, $N_{h}(\lambda)$ and
$\widehat{N}_h(\lambda)$ denote the number of $\lambda$-clusters for
$p$, $p_{h}$ and $\widehat{p}_{h}$, respectively. Notice that we do
not require $p$ to satisfy any smoothness properties. See Section
\ref{sec:cases} for the case case of smooth densities.
\begin{lemma}\label{lem:clust.numb}
Under conditions \textup{(C2)}
and for all $\varepsilon\in(0, \overline{\varepsilon})$ and $h \in(0,
\overline{h})$,
on the event $\mathcal{E}_{h,\varepsilon} = \{ \| \widehat{p}_h - p_h \|
_{\infty}< \varepsilon\}$,
\[
N_{h}(\lambda) = \widehat{N}_h(\lambda) = k.
\]
\end{lemma}

Condition (C3) is used to obtain establish rates of convergence for
the level set risk and the excess mass risk. It provides a way of
quantifying the clustering bias due to the use of the mollified
density $p_h$ as a function of the bandwidth $h$, locally in a
neighborhood of $\lambda$. In fact, if condition (C2) holds, then the
clustering bias is due to the sets $L_h(\lambda- \varepsilon) -
L(\lambda- \varepsilon)$, for $h \in(0,\overline{h})$ and $\varepsilon
\in
[0,\overline{\varepsilon})$.
\begin{lemma}\label{lem:xi.theta}
Under conditions \textup{(C2)} and \textup{(C3)}, for all $h \in
(0,\overline{h})$ and
$\varepsilon\in[0,\overline{\varepsilon})$ such that $\lambda-
\varepsilon
\geq0$,
%
%
\begin{equation}\label{eq:theta}
\mu\bigl( L_h(\lambda- \varepsilon) - L(\lambda- \varepsilon) \bigr)
\leq C_2 h^\theta,
\end{equation}
where
\[
\theta= \cases{
\displaystyle d - \max_i d_i + 1, &\quad if $\displaystyle\max_i d_i >
0$,\cr
d, &\quad otherwise,}
\]
and, for some positive constant $C_3$,
%
%
\begin{equation}\label{eq:xi}
P\bigl(L_h(\lambda- \varepsilon) - L(\lambda- \varepsilon)\bigr) \leq
C_3 h^\xi,
\end{equation}
where $\xi$ is either $\infty$ or $1$; in particular, $\xi= 1$ only
if $\max_i d_i = d$.
\end{lemma}

Condition (C3) is rather mild and depends only on dimension of the
support of~$P$. Indeed, it follows from the rectifiability property
(see the \hyperref[app]{Appendix} for details) that if, $S_i$ is a
component of the support of $P$ of dimension $d_i < d$, then $S_i$ has
box-counting dimension $d_i$
[see, e.g., \citet{ambrosiofuscopallara00}, Theorem 2.104]. This
implies that (C3) is satisfied for all $h$ small enough [see
also \citet{falconer03}]. For clusters $C_j$ belonging to
full-dimensional components of the
support of $P$, (C3) follows if the sets $\partial C_j(\lambda')$ have
box-counting dimension $d-1$ for all $\lambda' \in(\lambda-
\overline{\varepsilon}, \lambda]$ and if $\overline{h}$ is small
enough. In fact, under these additional assumptions, it is possible to
show that the bounds in Lemma \ref{lem:xi.theta} are sharp in the
sense that, for all $\lambda' \in(\lambda- \overline{\varepsilon},
\lambda]$ and $h \in(0, \overline{h})$, $\mu( L_h(\lambda') -
L(\lambda') ) = \Omega(h^\theta)$. In addition, provided that
$\overline{h}$ is smaller than the minimal inter-cluster distance
\[
{\min_{i \neq j} \inf_{x \in C_j y \in C_j}} \| x - y\|,
\]
we also obtain that $P(L_h(\lambda- \varepsilon) - L(\lambda-
\varepsilon
)) = \Omega(h^\xi)$, where $\xi$ can only be $1$ or~$\infty$.


Finally, we point out that the value of $\overline{h}$ depends on the
curvature of the components $\partial L(\lambda- \varepsilon)$ for all
$\varepsilon\in[0, \overline{\varepsilon})$, and on the minimal
inter-cluster distance. The smaller the condition numbers [see,
e.g., \citet{niyogismaleweinberger2008}] of these components,
and the
larger the inter-cluster distance, the larger $\overline{h}$.




Although the rates are not affected by the constants, in
practice, they can have a significant effect on the results, since
they may very well depend on $d$. This is especially true of $C_1$,
as illustrated in Example \ref{ex:gaussians} below.

\subsection{A refined analysis of condition \textup{(C1)}}\label{sec:C1}

We conclude this section with some comments on the parameter
$\gamma$ appearing in condition (C1), whose value affects in a
crucial way the consistency rates, with faster rates arising from
larger values of $\gamma$. If $S$ has dimension smaller than $d$,
then, clearly, $\gamma= \infty$, thus throughout this subsection we
assume that $P$ is a probability measure on $\mathbb{R}^d$ having
Lebesgue density $p$.

First, a
fairly general sufficient condition for assumption (C1) to hold with
$\gamma=1$ at $\lambda$ can be easily obtained using probabilistic
arguments as follows. Let $G$
denote the distribution of the random variable $Y=p(X)$ and suppose
$G$ has a Lebesgue density $g$ which is bounded away from $0$ and
infinity on $(\lambda-\overline{\varepsilon}, \lambda+
\overline{\varepsilon})$. Then, by the mean value theorem, for any nonnegative
$\varepsilon< \overline{\varepsilon}$,
\[
\mathbb{P}\bigl(\lambda- \varepsilon\leq p(X) \leq\lambda+ \varepsilon
\bigr) =
G\bigl(\{ y \dvtx y \in(\lambda+ \varepsilon,\lambda- \varepsilon) \}
\bigr)=
\varepsilon g(\lambda+ \eta)
\]
for some $\eta\in(- \varepsilon, \varepsilon)$. Thus, (C1) holds with
$\gamma= 1$ at $\lambda$. See also Example \ref{ex:gaussians}
below. A more refined result based on analytic conditions is given
next. Below $\mathcal{H}^{d-1}$ denotes the $(d-1)$-dimensional
Hausdorff measure in $\mathbb{R}^d$. See the \hyperref[app]{Appendix}
for the
definition of Hausdorff measure.
\begin{lemma}\label{lem:1}
Suppose that $P$ is a probability measure on $\mathbb{R}^d$ having
Lipschiz density $p$. Assume that, almost everywhere $\mu$, $\|
\nabla p(x) \| > 0$ and that $\mathcal{H}^{d-1}(\{ x \dvtx p(x) =
\lambda\}) < \infty$ for any $\lambda\in(0,\|
p\|_{\infty})$. Then \textup{(C1)} holds with $\gamma= 1$ for each
$\lambda\in(0,\| p\|_{\infty})$ except for a set of Lebesgue
measure $0$.
\end{lemma}

A further point of interest is to characterize the set of $\lambda$
values for which condition (C1) holds with $\gamma\neq1$. Clearly,
if $p$ has a jump discontinuity, then (C1) holds with $\gamma=
\infty$, for all values of $\lambda$ in some interval. On the other
hand, due to the previous result, if $ \| \nabla p \|$ is
bounded away from $0$ and $\infty$ in a neighborhood of
$p^{-1}(\lambda)$, then $\gamma= 1$. Thus, one could expect a value of
$\gamma$ different than $1$ when $\nabla p$ does not exist or when $\|
\nabla p\|$ is infinity or vanishes in $p^{-1}(\lambda)$. See the
example on page 7 in \citet{rigvert2006}, where (C1) holds with
$\gamma< 1$ if $q > d$ and $\gamma> 1$ if $q < d$, the former case
corresponding to $\| \nabla p(x_0) \| = 0$ and the latter to $\lim_{x
\rightarrow x_0} \| \nabla p(x) \| = \infty$. However, this would
seem to indicate that, if $p$ is sufficiently regular, the values of
$\lambda$ for which $\gamma\neq1$ form a negligible set of
$\mathbb{R}$. Lemma~\ref{lem:1} above already shows that this set has
Lebesgue measure zero if $p$ is Lipschitz with nonvanishing
gradient. Under stronger assumptions, it can be verified that this set
is in fact finite.
\begin{corollary}\label{cor:2}
Under the assumption of Lemma \ref{lem:1}, if $p$ is of class
$\mathcal{C}^1$ and has compact support, then the set of $\lambda$
such that \textup{(C1)} holds with $\gamma\neq1$ is finite.
\end{corollary}
\begin{example} \label{example::sharp}
Sharp clusters and lower-dimensional clusters. Suppose that
$p=\frac{dP }{d \mu} = \sum_{i=1}^m \pi_j p_j$
where
$p_i$ is a density with support on a compact, connected set $S_i$,
$\sum_i \pi_i =1$ and
$\min_i \pi_i >0$.
Moreover, suppose that
\[
{\min_{i \neq j} \inf_{x \in C_i, y \in C_j}}\| x - y \| > 0,
\]
where
$d(A,B) = {\inf_{x\in A,y\in B}} \Vert x-y\Vert$.
Finally, suppose that
\[
\min_j \inf_{x\in C_j} \pi_j p(x) \geq\lambda.
\]
We denote clusters of this type as sharp clusters. See \citet
{nowak2009},
for example.
It is easy to see that
(C1) and (\ref{eq:xi}) hold with $\gamma= \xi= \infty$.
A more general example in which one of the mixture component is
supported on a lower dimensional set is shown in Figure \ref{fig::sharp}.
%
%
\begin{figure}

\includegraphics{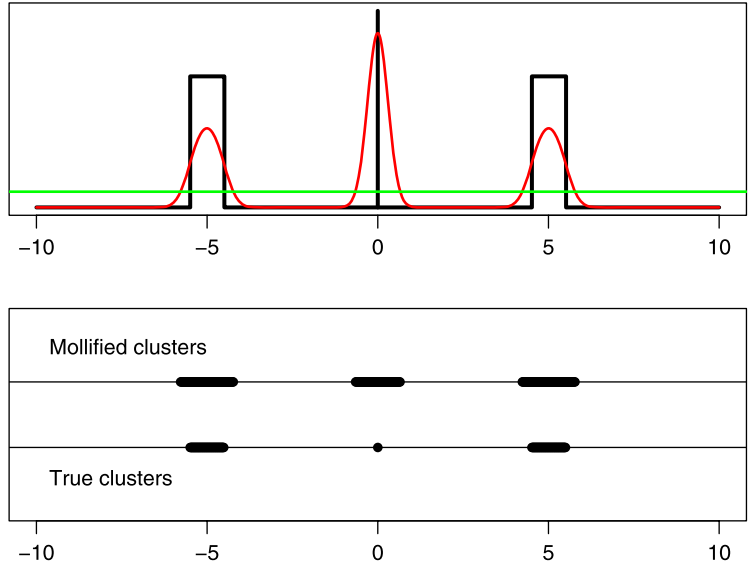}

\caption{Sharp clusters.
Top: the density of
$P = (1/3)\operatorname{Unif}(-5.5,-4.5)+ (1/3)\operatorname
{Unif}(4.5,5.5)+\break(1/3)\delta_0$
and the mollified density $p_h$ for $h=0.04$.
The point mass at 0 is indicated with a vertical bar.
Bottom: the true clusters and the
mollified clusters of $p_h$ with $\lambda= 0.04$.}
\label{fig::sharp}
\end{figure}
Here, the true distribution is
$P = (1/3)\operatorname{Unif}(-5.5,-4.5)+ (1/3)\operatorname
{Unif}(4.5,5.5)+(1/3)\delta_0$.
The geometric density and the mollified density based on $h=0.04$ are
shown in the top plot.
The point mass at $0$ is indicated with a vertical bar.
The bottom plot shows the true clusters and the mollified clusters
based on $p_h$ with $\lambda= 0.04$.
The clusters based on $p_h$ contain the true clusters and the difference
between them is a set of zero probability.
\end{example}
\begin{example}[(Normal distributions)] \label{ex:gaussians}
Suppose that $X \sim N_d(0,\Sigma)$, with $\Sigma$ positive definite.
Set $\sigma= |\Sigma|^{1/2}$. Then (C1) holds for any $0 \leq\lambda
\leq( \sigma(\sqrt{2 \pi} )^d )^{-1}$
with $\gamma=
1$ and $C_1 = C_d 2\sigma( \sqrt{2 \pi} )^d$, where the
constant $C_d$ depends on $d$ (and, of course, $\lambda$).
We prove the claim only for
$\lambda= \alpha(\sigma(\sqrt{2 \pi} )^d
)^{-1}$, where $\alpha\in(0,1)$.
Cases in
which $\alpha=1$ or $\alpha=0$ can be dealt with similarly.
Let $W \sim\chi^2_d$ and notice that $X^\top\Sigma^{-1} X
\stackrel{d}{=} W$. For all $\varepsilon> 0$ smaller than
%
%
\begin{equation}\label{eq:bound.e}
\min\biggl\{ \frac{\alpha}{\sigma( \sqrt{2 \pi}
)^d}, \frac{(1 - \alpha)}{\sigma( \sqrt{2 \pi}
)^d} \biggr\},
\end{equation}
simple algebra yields
\begin{eqnarray*}
&&P\bigl(|\phi_\sigma(X) - \lambda|< \varepsilon\bigr) \\
&&\qquad=
\mathbb{P} \biggl( 2 \log\frac{1}{\alpha- \varepsilon\sigma(
\sqrt{2 \pi} )^d } \leq W \leq
2 \log\frac{1}{\alpha+ \varepsilon\sigma( \sqrt{2 \pi}
)^d } \biggr)\\
&&\qquad= 2 \biggl( \log\frac{1}{\alpha- \varepsilon( \sigma\sqrt
{2 \pi} )^d} -
\log\frac{1}{\alpha+ \eta\sigma( \sqrt{2 \pi}
)^d} \biggr)
p_d \biggl( \log\frac{1}{\alpha+ \eta\sigma( \sqrt{2 \pi}
)^d} \biggr)
\end{eqnarray*}
for some $\eta\in(-\varepsilon,\varepsilon)$
where
$p_d$ denotes the density of a $\chi^2_d$ distribution and the second
equality holds in virtue of the mean value theorem.
By a first order Taylor expansion, for $\varepsilon\downarrow0$, the
first term on the right-hand side of the previous display can be
written as
\[
2 \varepsilon\sigma\bigl( \sqrt{2 \pi} \bigr)^d
\biggl( \frac{1}{\alpha- \varepsilon\sigma
( \sqrt{2 \pi} )^d} + \frac{1}{\alpha+ \varepsilon\sigma
( \sqrt{2 \pi} )^d} \biggr) + o(\varepsilon^2).
\]
Since $ ( \frac{1}{\alpha- \varepsilon\sigma
( \sqrt{2 \pi} )^d} +
\frac{1}{\alpha+ \varepsilon\sigma( \sqrt{2 \pi} )^d}
)
p_d ( \log\frac{1}{\alpha+ \eta\sigma( \sqrt{2 \pi}
)^d} ) \asymp1$
for any
$\varepsilon\geq0$ bounded by (\ref{eq:bound.e}), the claim is proved.
%
%
\begin{figure}

\includegraphics{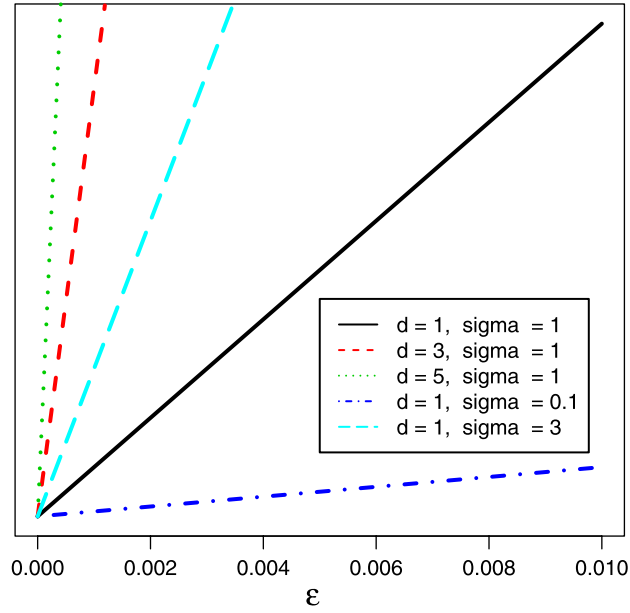}

\caption{Noise exponent for Gaussians.
Each curve shows $\mathbb{P}(|p(X) - \lambda| < \varepsilon)$ versus
$\varepsilon$ for $\alpha= 1/2$.
The plots are nearly linear since $\gamma=1$ in this case.}
\label{fig::normal_epsilon}
\end{figure}
See Figure \ref{fig::normal_epsilon}.
\end{example}

\section{Rates of convergence}
\label{section::rates}

In this section, we study the rates of convergence in the two distances
using deterministic bandwidths.
We defer the discussion of random (data driven)
bandwidths until Section
\ref{sec::bandwidth}.

\subsection{Preliminaries}\label{sec::preliminaries}
Before establishing consistency rates for the different risk measures
described above, we discuss some necessary preliminaries.

In our analysis, we require the event
%
%
\begin{equation}\label{eq:eventE}
\mathcal{E}_{h,\varepsilon} \equiv\{ \Vert\widehat{p}_h -
p_h\Vert_\infty\leq\varepsilon\},\qquad \varepsilon\in(0,\overline
{\varepsilon}), h \in(0,\overline{h}),
\end{equation}
to hold with high probability, for all $n$ large enough. In fact, some
control over $\mathcal{E}_{h,\varepsilon}$ provides a means of bounding
the clustering risks, as shown in the following result.
\begin{lemma}\label{lemma::handy2}
Let $\varepsilon\in(0, \overline{\varepsilon})$ and $h \in(0, \overline
{h})$ be such that the conditions \textup{(C1)} and \textup{(C2)} are
satisfied. Then, on
the event $\mathcal{E}_{h,\varepsilon}$,
\[
L(\lambda+ \varepsilon) \subseteq\widehat{L}_h(\lambda) \subseteq
L(\lambda+\varepsilon)\cup A \cup B,
\]
where
\[
A = L(\lambda- \varepsilon) - L(\lambda+ \varepsilon)
\]
and
\[
B = L_h(\lambda- \varepsilon) - L(\lambda- \varepsilon).
\]
Therefore, on $\mathcal{E}_{h,\varepsilon}$, under the additional
condition \textup{(C3)},
%
%
\begin{equation}\label{eq:firstrates}
P ( \widehat{L}_h(\lambda) \Delta L(\lambda) ) \leq
C_1\varepsilon^\gamma+ C_2 h^\xi.
\end{equation}
\end{lemma}

In order to bound $\mathbb{P}(\mathcal{E}^c_{h,\varepsilon})$, we study
the properties of the kernel estimator $\widehat{p}_h$.
We will impose the following condition on the kernel $K$.


\begin{enumerate}[(VC)]
\item[(VC)] The class of functions
\[
\mathcal{F} = \biggl\{ K \biggl( \frac{x - \cdot}{h} \biggr), x \in
\mathbb{R}^d, h > 0 \biggr\}
\]
satisfies, for some positive number $A$ and $v$
%
%
\begin{equation}\label{eq:VC}
\sup_P N\bigl(\mathcal{F}_h, L_2(P), \varepsilon\| F \|_{L_2(P)}\bigr)
\leq
\biggl( \frac{A}{\varepsilon} \biggr)^v,
\end{equation}
where $N(T,d, \varepsilon)$ denotes the $\varepsilon$-covering number
of the
metric space $(T,d)$, $F$ is the envelope function of
$\mathcal{F}$ and the supremum is taken over the set of all
probability measures on $\mathbb{R}^d$. The quantities $A$ and $v$ are
called the VC characteristics of $\mathcal{F}$.
\end{enumerate}
Assumption (VC) appears in \citet{Gine02},
\citet{einmhalmason2005}, and \citet
{ginekoltchinski06}. It holds
for a large class of kernels, including, for example, any compact
supported polynomial kernel and the Gaussian kernel. See
\citet{nolanpollard1987} and \citet{vdvwellner96} for
sufficient
conditions for (VC).

Using condition (VC), we can establish the following finite sample
bound for $\mathbb{P} ( \| \widehat{p}_h - p_h \|_{\infty} >
\varepsilon)$, which is obtained as a direct application of
results in \citet{Gine02}. 
%
\begin{proposition}[(Gine and Guillon)]
\label{prop::density}
Assume that the kernel satisfies the property \textup{(VC)} and that
%
%
\begin{equation}\label{eq:p}
\sup_{t \in\mathbb{R}^d} \sup_{h > 0} \int_{\mathbb{R}^d} K_h^2(t
- x) \,d P(x) < D < \infty.
\end{equation}
\begin{enumerate}
\item Let $h$ be fixed. Then, there exist constants $L > 0$ and $C >
0$, which depend only on
the VC characteristics of $K$, such that,\vspace*{1pt} for any $c_1
\geq C$ and $0 <
\varepsilon\leq
\frac{c_1 D}{\| K \|_{\infty}}$, there exists an $n_0 > 0$, which
depends on $\varepsilon$, $D$, $\| K \|_\infty$ and the VC
characteristics of $K$, such that, for all $n \geq n_0$,
%
%
\begin{equation}\label{eq:sup}\quad
\mathbb{P} \Bigl\{ {\sup_{x \in\mathbb{R}^d}}
|\widehat{p}_h(x) - p_h(x) | > 2 \varepsilon\Bigr\}
\leq
L\exp\biggl\{ - \frac{1}{L} \frac{\log(1 + c_1/(4L))}{c_1} \frac{n
h^d \varepsilon^2}{D} \biggr\}.
\end{equation}
\item Let $h_n \rightarrow0$ as $n \rightarrow\infty$ in such a way
that $\frac{n h_n^d}{|{\log h_n^d}|} \rightarrow\infty$. If $\{
\varepsilon_n \}$ is a sequence such that
%
%
\begin{equation}\label{eq:hn.en}
\varepsilon_n = \Omega\Biggl( \sqrt{\frac{\log r_n}{n h^d_n}} \Biggr),
\end{equation}
where $r_n = \Omega( h_n^{-d/2} )$, then, for all $n$
large enough, (\ref{eq:sup}) holds with $h$ and $\varepsilon$ replaced
by $h_n$ and $\varepsilon_n$, respectively. In particular, the term on
the right-hand side of (\ref{eq:sup}) vanishes at the rate $O
(r_n^{-1} )$.
\end{enumerate}
\end{proposition}

The above theorem imposes minimal assumptions on the kernel $K$ and,
more importantly, on the probability distribution $P$, whose density
is not required to be bounded or smooth, and, in fact, may not even
exist. Condition (\ref{eq:p}) is automatically satisfied by bounded
kernels. Finally, we remark that, for fixed $h$, setting $\varepsilon_n
= \sqrt{\frac{2 \log n}{h^d n C_K}}$ for an\vspace*{-2pt} appropriate
constant $C_K$
(depending on $K$), an application of the Borel--Cantelli lemma yields
that, as $n \to\infty$, $\| p_h - \widehat{p}_h \|_{\infty} =
O (
\sqrt{\frac{\log n}{n}} )$ almost everywhere
$P$. 

\subsection{Rates of convergence}

We now derive the converge rates for the clustering risks defined in
Section \ref{sec:risk}. Below, we will write $C_K$ for a constant whose
value depends only on the VC characteristic of the kernel $K$ and on
the constant $D$ appearing in (\ref{eq:p}).

We recall that Lemma \ref{lem:xi.theta} provides a way of controlling
the clustering bias due to the sets $L_h(\lambda- \varepsilon) -
L(\lambda- \varepsilon)$, uniformly over $\varepsilon<\overline
{\varepsilon}$
and $h < \overline{h}$. In fact, the parameters $\theta\in
\{1,\ldots,d \}$ and $\xi\in\{1,\infty\}$ will determine the rates
of consistency for the excess mass and the level set risk,
respectively. Specifically, higher values of the parameter $\theta$
which correspond to supports of lower dimension yield faster
convergence rates for the excess mass
risk. 
As for the level set risks, the case $\xi= \infty$ is the most
favorable, since it implies that the clustering bias has no effect on
the estimation of level sets and dimension independent rates are
possible. In particular, if $C_j$ has dimension smaller than $d$, then
$P(C_{\sigma(j)}^h - C_{j}) = 0$, so that $\xi= \infty$. More
generally, $\xi= \infty$ occurs when $L = S$. Overall our results
yield that, as expected, better rates for the clustering risk are
obtained for distributions supported on lower-dimensional sets.
\begin{theorem}[(Level set risk)]\label{thm::levelset}
Suppose that
\textup{(C1)}, \textup{(C2)}, \textup{(C3)} and \textup{(VC)} hold.
Then there exists a constant $C_L$ such that, for any $h \in(0,
\overline{h})$ and $ \varepsilon\in(0, \overline{\varepsilon})$,
%
%
\begin{equation}\label{eq:lev.set.risk}
R^L(p,\widehat{p}_h) \leq C_L ( \varepsilon^\gamma+ h^\xi+ e^{-
C_K n h^d \varepsilon^2} ).
\end{equation}
In particular, setting
\[
h_n = \biggl(\frac{ \log n}{n} \biggr)^{\gamma/(2\xi+ d\gamma)}
\quad\mbox{and} \quad\varepsilon_n = \sqrt{ \frac{\log n }{C_K n h_n^d}}
\]
we obtain
%
%
\begin{equation}
R^L(p,\widehat{p}_{h_n}) =
O \biggl( \max\biggl\{ \biggl( \frac{\log n}{n} \biggr)^{
{\gamma\xi}/({2\xi+ d\gamma})}, \frac{1}{n} \biggr\} \biggr).
\end{equation}
\end{theorem}

If $\gamma= \infty$, then either $S- L$ is empty or has zero Lebesgue
measure, or $S - L$ is a full dimensional set of positive Lebesgue
measure. The former cases, which correspond to $P$ having a
lower-dimensional support or to sharp clusters (see Example
\ref{example::sharp}), implies that $R^L =
O (\frac{1}{n} )$. Thus, we have dimension independent
rates for sharp clusters. In the latter case, $\xi=1$, so that $R^L$
is of order $O ( ( \frac{\log n}{n} )^{1/d}
)$.
When $\gamma< \infty$, then $\xi=1$ and the risk is of order $O
( ( \frac{\log n}{n } )^{{\gamma}/({2 + d
\gamma})} )$.

In practice, there are examples in between the sharp and nonsharp cases
for probability distributions with full-dimensional support.
For example, if there is a very small amount of mass just outside the cluster,
then, technically, $\xi=1$ and the rate will be slow for large $d$.
However, if this mass is very small then we expect
for finite samples that the
behavior of the risk will be close to the behavior observed in the
sharp case.
We could capture this idea mathematically by allowing $P$ to change
with $n$
and then allowing $\xi_n$ to vary with $n$ and take values between 1
and~$\infty$.
However, we shall not pursue the details here.

As an interesting corollary to Theorem \ref{thm::levelset}, we can
show that the
expected proportion of sample points that are incorrectly assigned as
clusters or noise vanishes at the same rate.
\begin{corollary}
\label{corollary::fI}
Let
$\widehat{f}_h = \frac{ |\widehat{I}_h|}{n}$,
where
\[
\widehat{I}_h = \bigl\{ i \dvtx\operatorname{sign}\bigl(\widehat
{p}_{h}(X_i) -
\lambda\bigr) \neq
\operatorname{sign}\bigl(p(X_i) - \lambda\bigr) \bigr\}.
\]
Then,
$\mathbb{E}(\widehat{f}_h) \leq C_L ( \varepsilon^\gamma+ h^\xi
+ e^{- C_K n h^d \varepsilon^2} )$.
\end{corollary}

We now turn to the excess mass risk.
\begin{theorem}[(Excess mass)]
\label{thm::excessmass}
Suppose that
\textup{(C1)}, \textup{(C2)}, \textup{(C3)} and \textup{(VC)} hold.
Then, there exists a constant $C_M$,
independent of $\varepsilon$ and $h$,
such that, for any $h \in(0, \overline{h})$ and $\varepsilon\in(0,
\overline{\varepsilon})$
with $\varepsilon< \lambda$,
\[
R^M(p,\widehat{p}_h) \leq C_M ( \varepsilon^{\gamma+1} +
h^\theta+ e^{-n C_K \varepsilon^2 h^d} ).
\]
Thus, setting
\[
h_n = \biggl(\frac{ \log n}{n} \biggr)^{({\gamma+ 1})/({2 \theta+
d(\gamma+ 1)})}
\quad\mbox{and}\quad \varepsilon_n = \sqrt{ \frac{\log n }{C_K n h_n^d}},
\]
we obtain
%
%
\begin{equation}
R^M(p,\widehat{p}_h) = O \biggl(
\biggl(\frac{\log n}{n} \biggr)^{{\theta(\gamma+ 1)}/({2 \theta
+ d(\gamma+ 1)})} \biggr).
\end{equation}
\end{theorem}

When $\gamma= \infty$ the excess mass risk $R^M$ is of order $O (
\frac{\log n}{n} )^{\theta/d}$. Thus, the higher~$\theta$,
that is, the smaller the dimension of the support of $P$, the faster the
rate of convergence. In particular, if $P$ is supported over a finite
set of points the~risk vanishes at the dimension independent rate $O
( \frac{\log n}{n} ) $. When $\gamma< \infty$, then
$\theta=1$ and the risk is of order $O ( ( \frac{\log n}{n
} )^{({\gamma+1})/({2 + d (\gamma+1)})})$.

\subsection{Some special cases}\label{sec:cases}

Here, we discuss some interesting special cases.

\subsubsection*{Fast rates for biased clusters}
In some cases, we might be content
with estimating the level set $L_h(\lambda)$, which is a biased
version of $L(\lambda)$. That is, the fringe $L_h(\lambda) -
L(\lambda)$ may not be of great practical concern and, in fact, it may
contain a very small amount of mass.
Indeed, we believe this is why clustering is often so successful in
high-dimensional problems.
Exact estimation of the level sets is not necessary in many practical problems.
In fact, by Lemma
\ref{lem:clust.numb}, conditions (C2) guarantees that $\widehat{L}$ will
include $L$ with high probability. Thus, for clustering purposes, one
may consider
some modifications of our risk functions.
First, suppose we only require that the estimated clusters
cover the true clusters.
That is, we say there is not error as long as $C_j\subset\widehat C_j$.
This suggests
the following modification of our risk functions:
\begin{itemize}
\item$\widetilde{R}^L(p,\widehat{p}_h) = \int_{ \{x \dvtx p(x)
\geq\lambda\} \cap\{x \dvtx\widehat{p}_h(x) < \lambda\}} dP(x)$,
\item$\widetilde{R}^M(p,\widehat{p}_h) = \mathcal{E}(L) - \mathbb
{E}(\widehat{L}_h \cap L)$.
\end{itemize}
Then we have the following result, which gives faster, dimension
independent rates.
The proof is similar to the proofs of the previous results and is omitted.
\begin{theorem}
\label{thm::biased}
Let $h \in(0, \overline{h})$ be fixed. Under \textup{(C1)}, \textup
{(C2)} and \textup{(VC)}, then
\[
\widetilde{R}^L(p,\widehat{p}_h)= O \biggl(\max\biggl\{
\biggl(\frac{\log n}{n} \biggr)^{\gamma/2}, \frac{1}{n} \biggr\}
\biggr)
\]
and
\[
\widetilde{R}^M(p,\widehat{p}_h) = O \biggl( \biggl(\frac{\log
n}{n} \biggr)^{({1+\gamma})/{2}} \biggr).
\]
\end{theorem}

Alternatively, one may be only interested in estimating the clusters
of the mollified density $p_h$, for any fixed $h \in
(0,\overline{h})$. Then, provided that $p_h$ is sufficiently smooth
(which is guaranteed by choosing a smooth kernel) and has finite
positive gradient for each point in the set $\partial L_h(\lambda)$,
the results in Section \ref{sec:C1} show that, for all $\varepsilon$
small enough,
\[
\mu\bigl( \{x \dvtx|p_h(x) - \lambda| < \varepsilon\} \bigr) \leq
\varepsilon.
\]
Thus, under assumptions (C2) and (VC), similar arguments to the ones
used in the proofs of Theorems \ref{thm::levelset} and \ref
{thm::excessmass} imply that
\[
R^L(p_h,\widehat{p}_h) = \int_{ \{x \dvtx p_h(x) \geq\lambda\}
\Delta\{x \dvtx\widehat{p}_h(x) \geq\lambda\}} dP(x) = O
\Biggl(\sqrt{\frac{\log n}{n}} \Biggr)
\]
and
\[
R^M(p_h,\widehat{p}_h) = \mathcal{E}(L_h) - \mathbb{E} (
\mathcal{E}(\widehat{L}_h) ) = O \biggl( \frac{\log n}{n}
\biggr).
\]
In either case, we get dimension independent rates.

\subsubsection*{The smooth full-dimensional case}
In the more specialized settings in which $P$ has full-dimensional
support and the Lebesgue density $p$ is smooth, better results are
possible. For example, using the same settings of
\citet{rigvert2006}, if $p$ is $\beta$-times H\"{o}lder
differentiable, then the bias conditions (C2) are superfluous, as
%
%
\begin{equation}\label{eq:p.bias}
\|
p_h - p \|_{\infty} \leq C h^\beta
\end{equation}
for some constant $C$ which
depends only on the kernel $K$. Choosing $h$ such that $C h^\beta<
\varepsilon$, on the event $\mathcal{E}_{h,\varepsilon}$, the triangle
inequality yields $\| \widehat{p}_h - p \|_\infty< 2 \varepsilon$. Thus,
for each $\varepsilon< \frac{\overline{\varepsilon}}{2}$ and each $h$ such
that $C h^\beta< \varepsilon$, on $\mathcal{E}_{h,\varepsilon}$,
instead of
(\ref{eq:firstrates}), one obtains
\[
P ( \widehat{L}_h(\lambda) \Delta L(\lambda) ) \leq C_1
2^\gamma\varepsilon^\gamma.
\]
Then, setting
$h_n = ( \log n/n )^{{1}/({2 \beta+ d})}$ and
$\varepsilon_n = \Omega( ( \log n/n))^{{\beta}/({2 \beta+ d})})$,
we see that
$R^L(p,\widehat{p}_h)$ is of order
$O(( \log n/n )^{{\gamma\beta}/({2 \beta+ d})})$, while
$R^M(p,\widehat{p}_h)$ is of order
$O((\log n /n )^{{(\gamma+1) \beta}/({2\beta+ d})})$.
These, are, up to an extra logarithmic factor, the minimax
rates established by \citet{rigvert2006}. In fact, under these
smoothness assumptions, and since the bias can be uniformly controlled
as in (\ref{eq:p.bias}), then, by
a combination of Fubini's theorem and of a peeling argument as in
\citet{audiberttsybakov2007} and \citet{rigvert2006}, the
exponential term $O ( e^{- C_K n h^d \varepsilon^2} )$
becomes redundant and rates without
the logarithmic term are possible.

\section{Choosing the bandwidth}
\label{sec::bandwidth}

In this section, we discuss two data-driven method for choosing the
bandwidth that adapts to the unknown parameters $\gamma$ and~$\theta$.
Before we explain the details, we point out that $L_2$
cross-validation is not appropriate for this problem.
In fact, we are allowing for the case where $P$ may have atoms,
in which case it is well known that cross-validation chooses $h=0$.

\subsection{Excess mass}

We propose choosing $h$
by splitting the data
and maximizing an empirical estimate of the
excess mass functional.
Polonik (1995) used this approach to
choose a level set from among a fixed class
$\mathcal{L}$ of level sets of finite VC dimension.
Here, we are choosing a bandwidth, or, in other words,
we are choosing a level set from
a random class of level sets
$\mathcal{L}= \{ \{x \dvtx\widehat{p}_h(x) \geq\lambda\}\dvtx h > 0\}
$ depending on the observed sample $X$.
The steps are in Table
\ref{fig::bandwidth}.

%
\begin{table}[b]
\caption{Selecting the bandwidth using the excess mass risk}
\label{fig::bandwidth}
\rule{\textwidth}{0.5pt}
\vspace*{-14pt}
\begin{enumerate}
\item Split the data into two halves which we denote by
$X=(X_1,\ldots,X_n)$ and
$Z=(Z_1,\ldots, Z_n)$.
\item {Let $\mathcal{H}$ be a finite set of bandwidths.
Using $X$, construct kernel density estimators
$\{\widehat{p}_h\dvtx h \in\mathcal{H}\}$.
Let $L_h = \{x\dvtx\widehat{p}_h(x) \geq\lambda\}$.}
\item {Using $Z$, estimate the excess mass functional}
\[
\widehat{\mathcal{E}}(h) = \frac{1}{n}\sum_{i=1}^n I(Z_i \in L_h) -
\lambda\mu(L_h).
\]
\item {Let $\widehat{h}$ be the maximizer of $\widehat
{\mathcal{E}}(h)$
and set $\widehat{L} = L_{\widehat{h}}$.}\vspace*{-8pt}
\end{enumerate}
\rule{\textwidth}{0.5pt}
%
\end{table}

To implement the method, we need to compute
$\mu(L_h)$.
In practice,
$\mu(L_h)$ can be approximated by
\[
\frac{1}{M}\sum_{i=1}^M\frac{I(\widehat{p}_h(U_i) \geq\lambda)
}{g(U_i)},
\]
where $U_1,\ldots, U_M$ is a sample from a convenient density $g$.
In particular, one can choose $g=\widehat{p}_H$ for some large
bandwidth $H$.
Choosing $M \approx n^2$ ensures that the extra error of this
importance sampling estimator
is $O(1/n)$ which is negligible.
We ignore this error in what follows.

Technically, the method only applies for $\lambda>0$,
at least in terms of the theory that we derive.
In practice, it can be used for $\lambda=0$.
In this case, $\widehat{\mathcal{E}}(h)$ becomes 1 when $h$ is large.
We then take $\widehat{h}$ to be the smallest $h$ for which
$\widehat{\mathcal{E}}(h)=1$.

Below we use the notation
$\mathcal{E}_X(\cdot)$ instead of $\mathcal{E}(\cdot)$ to indicate
that the excess mass functional (\ref{eq:excess.mass}) is
evaluated at a random set depending on the training set
$X$ and, therefore, is itself random. Accordingly, with some abuse of
notation, for any \mbox{$h>0$}, we will
write $\mathcal{E}_X(h) = \mathcal{E}(L_h)$, with $L_h$ the
$\lambda$-level set of $\widehat{p}_h$. Below $\mathcal{H}$ is a
countable dense subset of $[0,\overline{h}]$.
The next result is closely related to Theorem 7.1 of \citet{gyorfi2002}.
\begin{theorem}\label{thm::split}
Let
$h_* = \arg\max_{h\in\mathcal{H}} \mathcal{E}_X(h)$.
For any $\delta>0$,
%
%
\begin{equation}\label{eq::expect}
\mathbb{E}(\mathcal{E}_X(h_*)) - \mathbb{E}(\mathcal{E}_X(\widehat
{h})) \leq
d(\delta,\kappa) \frac{1+\log2}{n},
\end{equation}
where the\vspace*{1pt} expectation is with respect to the
joint distribution of the training and test set, $d(\delta,\kappa) =
\frac{2}{\kappa} \delta(1 + \delta) (
16 \gamma^2 + \delta(7 + 16 \gamma^2) )$,
with $\kappa= 2 + \lambda\mu(S + B(0,\overline{h}))$ and $\gamma^2
= \frac{7}{4} ( e^{ 4/7} - 1 )$.
\end{theorem}

Now we construct a grid $\mathcal{H}_n$ of size depending on $n$ that
is guaranteed
to ensure that optimizing over $\mathcal{H}_n$ implies we are adapting
over $\gamma$ and $\theta$.
\begin{theorem}
\label{thm::adapt}
Suppose \textup{(C1)} and \textup{(C2)} hold.
Let
\[
\delta_n(\theta) = \frac{a_n^{\theta/d}}{2 A_n(\theta)},
\]
where
$a_n=(\log n/n)$ and
\[
A_n(\theta) = \frac{2 |{\log a_n}| a_n^{\theta/(2\theta+d)}\theta
^2}{(2\theta+d)^2}.
\]
Let
$G_n(\theta) = \{\gamma_1(\theta),\ldots, \gamma_{N(\theta
)}(\theta)\}$
where $\gamma_j(\theta) = (j-1)\delta_n(\theta)$ and
$N(\theta)$ is the smallest integer less than or equal to
$\Upsilon_n(\theta)/\delta_n(\theta)$,
\[
\Upsilon_n(\theta) = \frac{2\theta^2}{d^2 W_n} - \frac{2\theta
}{d} -1
\]
and
\[
W_n = \frac{\log2}{\log n - \log\log n}.
\]
Let
\[
\mathcal{H}_n = \bigl\{ h_n(\gamma,\theta)\dvtx\theta\in\{1,\ldots, d\},
\gamma\in G_n(\theta)\bigr\},
\]
where
$h_n(\gamma,\theta) = a_n^{ (\gamma+1)/(2\theta+ d(\gamma+1))}$.
Let $\widehat L$ be obtained by minimizing $\widehat{\mathcal{E}}(h)$
for $h\in\mathcal{H}_n$.
Then
\[
\mathcal{E}(L) - \mathbb{E}(\mathcal{E}(\widehat{L})) \leq
O \biggl(\frac{\log n}{n} \biggr)^{{\theta(\gamma+1)}/({2\theta+
d(\gamma+1)})}.
\]
\end{theorem}

The latter theorem shows that our cross-validation
methods gives a completely data-driven method for choosing the
bandwidth that preserves the rate. Notice, in particular, that adapting
to the parameter $\theta$ is equivalent to adapting to the unknown
dimension of the support of $P$.
This makes it possible to use our method in practical problems
as long as the sample size is large.
For small sample sizes, data splitting might lead to highly variable results
in which case our bandwidth selection method might not
work well.
An alternative is to split the data many times and combines the estimates
over multiple splits.

When $\mu(L)=0$, we have that
$h_* =0$.
The above theorems are still valid in this case.
Thus, the case where $P$ is atomic is included while it is ruled out
for $L_2$ cross-validation.

\subsection{Stability}
\label{section::stability}

Another method for selecting the bandwidth is to
choose the value for $h$ that produces stable clusters, in a sense
defined below.
The use of stability has gained much popularity in clustering;
see
\citet{benhur2002} and \citet{lange2004},
for example.
In the context of $k$-means clustering and related methods,
\citet{bendavid2006}
showed that
minimizing instability leads to
poor clustering.
Here, we investigate the use of stability for density clustering.

Suppose, for simplicity, that the sample size is a multiple of 3.
That is, the sample size is $3n$ say. Now randomly split the data
into three vectors of size $n$, denoted by
$X=(X_1,\ldots,X_n)$,
$Y=(Y_1,\ldots,Y_n)$ and
$Z=(Z_1,\ldots,Z_n)$. (In practice, we split the data into three
approximately equal subsets.)

We define the \textit{instability function} as the random function $\Xi
\dvtx[0, \infty) \mapsto[0,1]$ given by
%
%
\begin{equation}
\Xi(h) \equiv\rho(\widehat{p}_h,\widehat{q}_h,\widehat{P}_Z)
=\int_{ \{ x \dvtx\widehat{p}_h(x) \geq\lambda\} \Delta\{ x
\dvtx\widehat{q}_h(x) \geq\lambda\}} d \widehat{P}_Z(x),
\end{equation}
where $\widehat{p}_h$ is constructed from $X$,
$\widehat{q}_h$ is constructed from $Y$ and
$\widehat{P}_Z$ is the empirical distribution based on $Z$.

Rather than studying stability in generality,
we consider a special case
involving the following extra conditions.

\begin{enumerate}
\item\textit{Sharp clusters.}
Assume that
$P = \sum_{j=1}^m \pi_j P_j$
where $\sum_i \pi_j =1$,
and $P_j$ is uniform on the compact set $S_j$ of full dimension $d$.
Thus,
$p(z) = \sum_j \Delta_j I(z\in S_j)$
where
$\Delta_j = \pi_j/\mu(S_j)$.
Let $\underline{\Delta} = \min_j \Delta_j> 0$ and
let $\overline{\Delta} = \max_j \Delta_j$.
\item\textit{Spherical Kernel.}
We use a spherical kernel
so that
\[
\widehat{p}_h(z) = \frac{1}{n h^d} \sum_{i=1}^n
\frac{I(\Vert z-X_i\Vert\leq h)}{v_d} =
\frac{\widehat{P}(B(x,h))}{h^d v_d},
\]
where
$v_d=\pi^{d/2}/\Gamma(d/2+1)$ denotes the volume of the unit ball
and
$\widehat{P}$ is the empirical measure.
\item\textit{The support of $P$ is a standard set.} Letting $S = \cup
_{j=1}^m S_j$, we assume that there exists a $\delta\in(0,1)$ such that
\[
\mu\bigl(B(z,h)\cap L\bigr) \geq\delta\mu(B(z,h))\qquad
\mbox{for all }z\in S \mbox{ and all } h < \operatorname{diam}(S),
\]
where $\operatorname{diam}(S) = \sup_{(x,y) \subset S } \| x - y \| $
indicates the diameter of the set $S$. This property appears in a
natural way in set estimation problems; see, for example, \citet
{cuevasfraiman1997}.
\item\textit{Choice of $\lambda$.} We take $\lambda= 0$, so that $L = S$.
\end{enumerate}

Under these settings, the graph $\Xi(h)$ is typically unimodal
with $\Xi(0) = \Xi(\infty) = 0$.
Hence, it makes no sense to minimize $\Xi$.
Instead, we will fix a constant $\alpha\in(0,1)$
and choose
%
%
\begin{equation}\label{eq::alpha}
\widehat{h} = \inf\Bigl\{h\dvtx\sup_{t>h}\Xi(t) \leq\alpha
\Bigr\}.
\end{equation}
\begin{theorem}
\label{thm::stability}
Let $h_* = \operatorname{diam}(L)$. Under conditions 1--4:
\begin{enumerate}
\item$\Xi(0) = 0$ and $\Xi(h) = 0$, for all $h \geq h_*$;
\item$\sup_{0 < h < h_*} \mathbb{E}(\Xi(h)) \leq1/2$;
\item as $h \rightarrow0$, $\mathbb{E}(\Xi(h)) \asymp h^d$;
\item for each $h \in(0, h_*)$,
\[
D_3 (h_*-h)^{d(n+1)} D_4^n \leq\mathbb{E}(\Xi(h)) \leq2 D_1
(h_*-h)^{n+1} D_2^n,
\]
where
\begin{eqnarray*}
D_1 &=& \frac{\pi^{d/2} h_*^{d-1}}{2^d \Gamma((d/2)+1)},\qquad
D_2 = \frac{\pi^{d/2} h_*^{d-1}}{\Gamma((d/2)+1)},\\
D_3 &=& \frac{\delta\underline{\Delta}\pi^{d/2} }{\Gamma
((d/2)+1)},\qquad
D_4 = \frac{\underline{\Delta} \delta\pi^{d/2}}{\Gamma(d/2+1)}.
\end{eqnarray*}
\end{enumerate}
\end{theorem}

To see the implication of Theorem \ref{thm::stability}, we
proceed as follows. Consider a grid of values $\mathcal{H} \subset
(0,h_*)$ of cardinality $n^\beta$, for some $0 <
\beta< 1$.
By Hoeffding's inequality, with probability at least $1 -
\frac{1}{n}$, we have that
\[
\sup_{h \in\mathcal{H}} | \Xi(h) - \mathbb{E}(\Xi(h)) | \leq w_n
\equiv\sqrt{\frac{2 \log(2 n)(1 - \beta) }{n}}.
\]
Replacing $\mathbb{E}(\Xi(h))$ by $\Xi(h) + w_n$ and $\Xi(h) - w_n$ in
the upper and lower bounds of part 4 of Theorem
\ref{thm::stability}, respectively, setting them both equal to
$\alpha$ and then finally solving for $h$, we conclude that the selected
$\widehat{h}$ is upper bounded by
\[
h_* - \biggl( \frac{\alpha- w_n}{2D_1} \biggr)^{1/(n+1)}
D_2^{ - {n}/({n+1})}
\]
and lower bounded by
\[
h_* - \biggl( \frac{\alpha+ w_n}{D_3} \biggr)^{1/(d(n+1))}
D_4^{ - {n}/({d(n+1)})}
\]
with probability larger than $1 - \frac{1}{n}$. Thus, as $n \to
\infty$, the resulting bandwidth does not tend to $0$. Hence, the
stability based method leads to bandwidths that are quite different than
the method in the previous section. Our explanation for this finding
is that the stability criterion is essentially aimed at reducing the
variability of the clustering solution, but it is virtually unaffected by
the bias caused by large bandwidths.

In the analysis above, we assumed for simplicity that $\lambda=0$.
When $\lambda>0$,
the instability $\Xi(h)$ can have some large peaks for very large $h$.
This occurs when $h$ is large enough so that
some mode of $p_h(x)$ is close to $\lambda$.
Choosing $h$ according to (\ref{eq::alpha}) will then lead to
serious oversmoothing.
Instead, we can choose $\widehat{h}$ as follows.
Let $h_0 =\arg\max_h \Xi(h)$ and define
%
%
\begin{equation}\label{eq::better}
\widehat{h} = \inf\{ h\dvtx h \geq h_0, \Xi(h)\leq\alpha
\}.
\end{equation}
We will revisit this issue in Section \ref{section::examples}.
A theoretical analysis of this modified procedure is tedious and,
in the interest of space, we shall not pursue it here.

\section{Approximating the clusters}
\label{sec::friend}

Lemma \ref{lem:clust.numb} shows that, under mild, conditions and when
the sample size is large enough, $N(\lambda) = \widehat{N}_h(\lambda)$
uniformly over $h \in(0,\overline{h})$ with high
probability. However, computing the number of connected components of
$\widehat{L}_h(\lambda)$ exactly is computationally difficult, especially
if $d$ is large.
In this section, we study a graph-based algorithm for finding the connected
components of $\widehat{L}_h$ and for estimating the number of
$\lambda$-clusters $N(\lambda)$ that is based on the $\rho$-nearest
neighborhood graph of $\{ X_i \dvtx\widehat{p}_h(X_i) \geq\lambda\}
$ that is fast and easy to implement.

The idea using the union of balls of radius $\rho$ centered at the
sample points to recover certain properties of the support of a
probability distribution is well understood. For instance, \citet
{devroywse80} and \citet{tsybakovkorostelev93} use it as a
simple yet effective estimator of the support, while \citet
{niyogismaleweinberger2008} show how it can be utilized for
identifying certain homology features of the support.

In particular, \citet{cuevas2000} and \citet{cadre2007}
propose to
combine a kernel density estimation with a single-linkage graph
algorithm to estimate the number of $\lambda$-clusters [see also
\citet{janghendry07}, for an application to large databases].
Our results offer similar guarantees but hold under more general
settings.

The algorithm proceeds as follows. For
some $h \in(0,\overline{h})$ and a given $\lambda\geq0$:
\begin{enumerate}
\item compute the kernel density estimate $\widehat{p}_h$;
\item compute the $\rho$-nearest neighborhood graph of $\{ X_i \dvtx
\widehat{p}_h(X_i) \geq\lambda\}$, that is the graph $\mathcal{G}_{h,n}$
on $\{ X_i \dvtx\widehat{p}_h(X_i) \geq\lambda\}$ where there is
an edge
between any two nodes if and only if they both belong to a ball of
radius $\rho$;
\item compute the connected components of $\mathcal{G}_{h,n}$ using a
depth-first search.
\end{enumerate}
The computational complexity of the last step is linear in the number
of nodes and the number of edges of $\mathcal{G}_{h,n}$ [see,
e.g., \citet{algo02}], which are both random.

We will show that, if $\rho$ is chosen appropriately, then, with high
probability as $n \to\infty$:
\begin{enumerate}
\item the number of connected components of $\mathcal{G}_{h,n}$,
$\widehat{N}^G_h(\lambda)$, matches the number of true clusters,
$N(\lambda) = k$;
\item there exists a permutation of $\{1, \ldots,k \}$ such that, for
each $j$ and $j'$,
%
%
\begin{equation}\label{eq:supp}\quad
C_j^h \subseteq\bigcup_{x \in\mathcal{C}_{\sigma(j)}}B(x,\rho)
\quad\mbox{and}\quad \biggl( \bigcup_{x \in\mathcal{C}_{\sigma
(j)}}B(x,\rho) \biggr) \cap
\biggl( \bigcup_{x \in\mathcal{C}_{\sigma(j')}}B(x,\rho) \biggr)
= \varnothing,
\end{equation}
where $\mathcal{C}_1,\ldots,\mathcal{C}_k$ are the connected
components of $\mathcal{G}_{h,n}$.
\end{enumerate}

We will assume the following regularity condition on the densities
$p_h$, which is satisfied if the kernel $K$ is of class
$\mathcal{C}^1$ and $P$ is not flat in a neighborhood of $\lambda$:
\begin{enumerate}[(G)]
\item[(G)] there exist constants $\varepsilon_1 > 0$ and $C_g > 0$ such
that for each $h \in(0,\overline{h})$, $p_h$ is of class
$\mathcal{C}^1$ on $\{ x \dvtx| p_{h}(x) -\lambda| < \varepsilon_1
\}$ and
%
%
\begin{equation}\label{eq:grad}
\inf_{h \in(0,\overline{h})} \inf_{x \in\{| p_h(x) - \lambda| <
\varepsilon_1 \}} \| \nabla p_h(x) \| > C_g.
\end{equation}
\end{enumerate}
Let $\delta_h = \min_{i \neq j} \inf_{x \in C^h_i,y \in C^h_j}\| x -
y\|$ and set
$\delta= \inf_{h \in(0,\overline{h})} \delta_h$. Notice that, under
(C2)(b), $\delta> 0$. Finally, let $\mathcal{O}_{h,n}$ denote the
event in (\ref{eq:supp}), which clearly implies the event $\{
\widehat{N}^G_h(\lambda) = k \}$.
\begin{theorem}
\label{thm::friends}
Assume conditions \textup{(G)} and \textup{(C2)} and let $d^* =
\operatorname{dim}(L)$. Assume further that there exists a constant
$\overline{C}$ such that, for every $r \leq\delta/2$ and for
$P$-almost all $x \in S \cap L$,
%
%
\begin{equation}\label{eq:standard}
P(B(x,r)) > \overline{C} r^{d_i},
\end{equation}
where $d_i = \operatorname{dim}(S_i)$, with $x \in S_i$. Then there exists
positive constants $\overline{\rho}$ and $\overline{M}$, depending
on $d^*$ and $L$ such that, for every $\rho< \min\{ \delta/2,
\overline{\rho}\}$, there exists a number $\varepsilon(\rho)$ such
that, for any $\varepsilon< \eta(\rho)$,
\[
\mathbb{P} ( \mathcal{O}_{h,n}^c ) \leq
\mathbb{P}(\mathcal{E}^c_{h,\varepsilon}) + \overline{M} \rho^{-d^*}
e^{- \overline{C} n \rho^{d^*}},
\]
uniformly in $ h \in(0, \overline{h})$.
\end{theorem}

The previous result deserves few comments. First, the constants
$\overline{\rho}$, $\overline{M}$ and $\overline{C}$ depend on $d^*$.
Second, assumption (\ref{eq:standard}) is a natural generalization
to lower-dimensional sets of the standardness assumption used, for
example, in
\citet{cuevasfraiman1997}. It is clearly true for components $P_i$
of full-dimensional support that are absolutely continuous with
respect to the Lebesgue measure. Finally, in view of Lemma
\ref{lem:suff.cond} [and, specifically, of the way
$\varepsilon(\rho,\tau)$ is defined], Theorem \ref{thm::friends} holds
for sequences $\{ \varepsilon_n \}$, $\{ h_n \}$ and $\{ \rho_n \}$ such
that:
\begin{enumerate}
\item$\varepsilon_n = o(1)$,
\item$\sup_n h_n \leq\overline{h}$;
\item$\sup_n \rho_n < \min\{ \delta/2, \overline{\rho}_d \}$ and
$\varepsilon_n = o(\rho_n)$.
\end{enumerate}
In particular,\vspace*{-2pt} if $h_n = o(1)$, then, following Proposition \ref
{prop::density}, the term $\mathbb{P}(\mathcal{E}^c_{h_n,\varepsilon
_n})$ vanishes if $\frac{n h_n^d}{|{\log h_n^d}|} \rightarrow\infty$.
Interestingly enough, condition (C1) does not play a direct role in
Theorem \ref{thm::friends}.

We now consider a bootstrap extension of the previous algorithm, as
suggested in \citet{cuevas2000}. For any $h$, let $X^* =
(X^*_1,\ldots,X^*_N)$, denote a bootstrap sample from $\widehat{p}_{h}$
conditionally on $\{ \widehat{p}_{h} \geq\lambda\}$ and let
$\mathcal{G}^*_{n,h}$ denote the $\rho$-neighborhood graph with node
set $X^*$. Finally, let $\mathcal{O}^*_{h,n}$ be the event given in
(\ref{eq:supp}), except that
$\mathcal{C}_1,\ldots,\mathcal{C}_k$ are now the connected components
of $\mathcal{G}^*_{h,n}$.
\begin{theorem}
\label{thm::bootstrap}
Assume conditions \textup{(C2)} and
\textup{(G)}. Suppose that there exist positive constants $\overline{C}$ and
$\overline{\rho}$ such that
%
%
\begin{equation}\label{eq:cond.boot}
\inf_{h \in(0,\overline{h})} \int_{A_h \cap L_h(\lambda) } p_h
\,d\mu>\overline{C} \rho^{d}
\end{equation}
for any ball $A_h$ of radius $\rho< \overline{\rho}$ and center in
$L_h(\lambda)$. 
Then, for any $\rho\leq\min\{ \delta/2, \overline{\rho} \}$,
there exists a positive number $\varepsilon(\rho)$ such that, for each
$\varepsilon< \varepsilon(\rho)$,
\[
\mathbb{P} ( (\mathcal{O}^*_{h,n})^c ) \leq
\mathbb{P}(\mathcal{E}^c_{h,\varepsilon}) + \overline{M} \rho^{-d}
e^{- C N \rho^{d} },
\]
uniformly in $ h \in(0, \overline{h})$, where $\overline{M}$ and $C$
are positive constants independent of $h$ and $\rho$.
\end{theorem}

The constants $C$, $\overline{C}$, $\overline{\rho}$, and
$\overline{M}$ depend on both $d$ and $S \oplus B(0,\overline{h})$.
In our
settings, condition (\ref{eq:cond.boot}) clearly holds if $P$ has
full-dimensional support. More generally, it can be shown that
conditions (G) and (\ref{eq:standard}) imply (\ref{eq:cond.boot}).

Just like with Theorem \ref{thm::friends},
using Lemma \ref{lem:suff.cond}, it can be verified that the theorem
holds if one consider sequences of parameters depending on the sample
size such that $\varepsilon_n = o(1)$, $\varepsilon_n = o(\rho_n)$,
$\sup_n
\rho_n < \max\{ \delta/2,\overline{\rho}\}$ and $\sup_n h_n <
\overline{h}$, provided that the conditions of Proposition
\ref{prop::density} are met.

Despite the similar form for the error bounds of Theorems \ref
{thm::friends} and \ref{thm::bootstrap}, there are some marked
differences. In fact,
in Theorem \ref{thm::friends} the performance of the algorithm
depends directly on the sample size $n$ and, in particular, on the
actual dimension $d^* \leq d$ of the support of $P$, with smaller
values of $d^*$ yielding better guarantees. In contrast, besides
$n$, the performance of the algorithm based on the bootstrap sample
depends on the ambient dimension $d$, regardless of $d^*$, and on
the bootstrap sample size $N$. By choosing $N$ very large, the
expression $\mathbb{P}(\mathcal{E}^c_{h,\varepsilon})$ becomes the
leading term in the upper bound of the probability of the event $
(\mathcal{O}^*_{h,n})^c$.

\section{Examples}
\label{section::examples}

In this section, we consider a few examples to illustrate the methods.

\subsection{A one dimensional example}

In Section \ref{section::stability},
we pointed out that when $\lambda>0$ and large,
it is safer to use the modified rule
$\widehat{h} = \inf\{ h\dvtx h \geq h_0, \Xi(h)\leq\alpha\}$
where $h_0 =\arg\max_h \Xi(h)$,
in place of the original rule
$\widehat{h} = \inf\{h\dvtx\sup_{t>h}\Xi(t) \leq\alpha\}$.
We illustrate this with a simple one-dimensional example.

%
\begin{figure}

\includegraphics{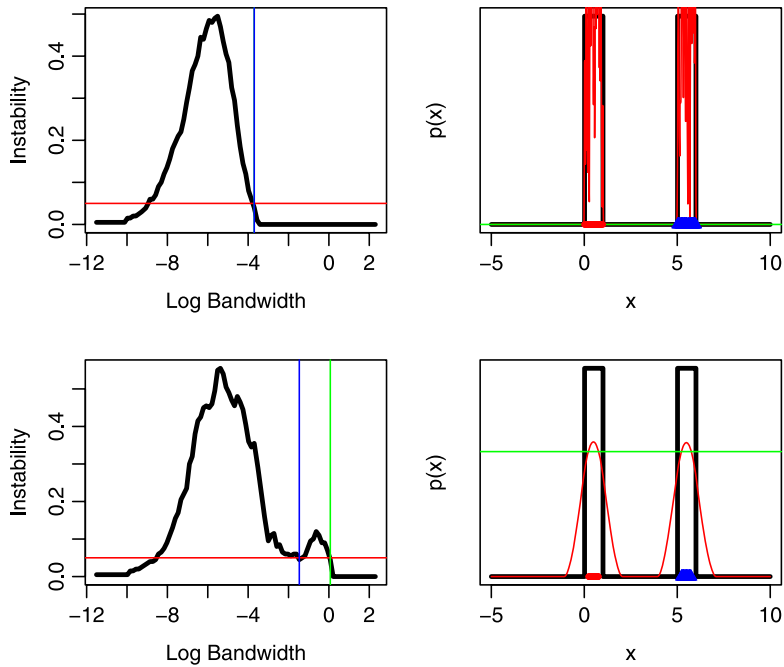}

\caption{The left plots show the instability as a function of log bandwidth.
The horizontal line shows $\alpha= 0.05$.
The right plots show the true density and the kernel density estimator
based on the selected bandwidth $h$.
In the top plots, $\lambda=0$.
In the bottom plots, $\lambda=0.3$.}
\label{fig::OneDim1}
\end{figure}

Figure \ref{fig::OneDim1}
shows an example based on
$n=200$ points from the density $p$ that
is uniform on
$[0,1]\cup[5,6]$.
When $\lambda=0$ (top),
the original rule
works fine.
(We use $\alpha$ = 0.05.)
The selected bandwidth is small leading to the very wiggly density
estimator in the top right plot.
However, this estimator correctly estimates the level set and the clusters.
In the bottom, we have $\lambda=0.3$.
When $h$ is large, there is a blip
in the instability curve
corresponding to the fact that the modes
of $p_h(x)$ are close to $\lambda$.
The original rule corresponds to the second vertical line in the
bottom left plot.
The resulting density estimator shown in the bottom right plot is
oversmoothed and leads to no points
being in the set $\widehat{p}_h \geq\lambda$.
The modified rule corresponds to the first vertical line
in the bottom left plot.
This bandwidth works fine.

%
\begin{figure}

\includegraphics{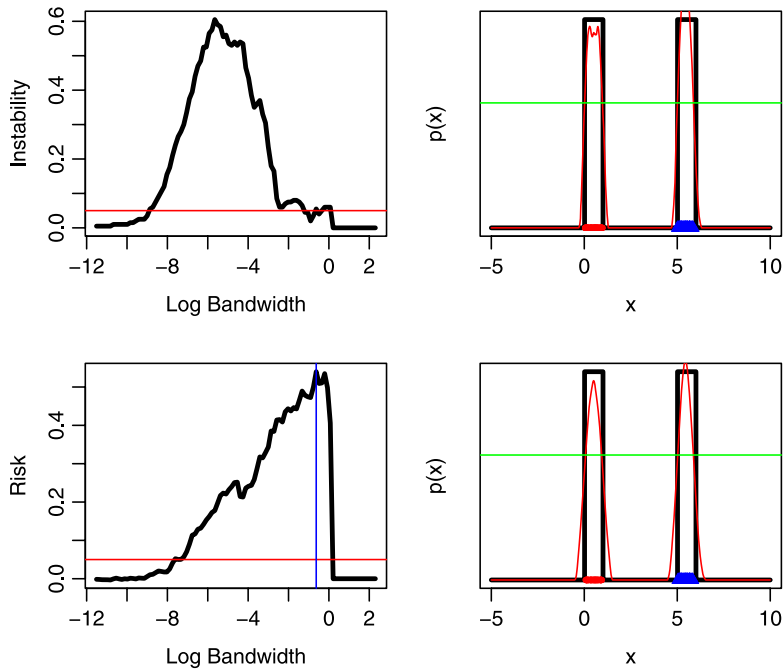}

\caption{The top left plot shows the instability as a function of log
bandwidth.
The top right plot shows the true density and the kernel density estimator
based on the selected bandwidth $h$ using the modified rule.
The bottom left plot shows the estimated excess mass risk as a function
of log bandwidth.
The top right plot shows the true density and the kernel density estimator
based on the selected bandwidth $h$ obtained by maximizing the excess mass.
In both bottom plots, $\lambda=0.3$.
Both methods recover the level set and the clusters.}
\label{fig::OneDim2}
\end{figure}

Figure \ref{fig::OneDim2}
compares the instability method (top)
with the
excess mass method (bottom).
Both methods recover the level set and the clusters.
We took $\lambda= 0.3$ in both cases.
Because $\lambda$ is very large,
the excess mass becomes undefined for large $h$
since $p_h(x) < \lambda$ for all $x$, which
we denoted by setting the risk to 0 in the bottom left plot.

\subsection{Fuzzy stick with spiral}

\begin{figure}

\includegraphics{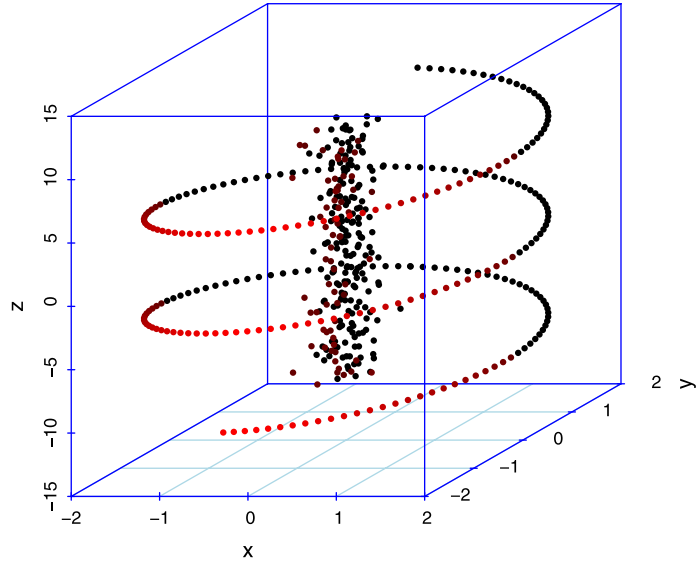}

\caption{500 data points from a fuzzy stick plus a spiral.}
\label{fig::fuzzy1}
\end{figure}

\begin{figure}[b]

\includegraphics{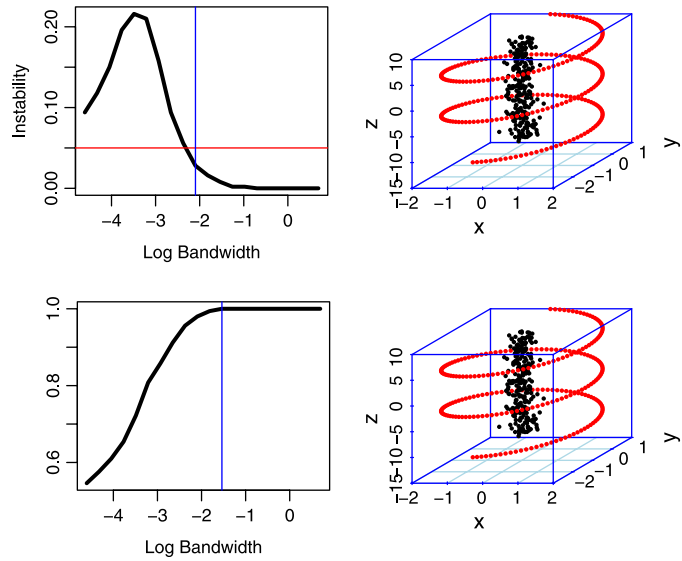}

\caption{Clusters obtained from instability (top)
and excess mass (bottom).}
\label{fig::fuzzy2}
\end{figure}

Figure \ref{fig::fuzzy1} shows data from
a fuzzy stick with a spiral.
The stick has noise while the spiral is
supported on a lower-dimensional curve.
Figure \ref{fig::fuzzy2} shows the clusterings from
the instability method and the excess risk method
with $\lambda=0$.
Both recover the clusters perfectly.
Note that the excess risk is necessarily
equal to 1 for large $h$.
In this case, we take $\widehat{h}$
to be the smallest $h$ of all bandwidths that
maximize the excess mass.
We see that both methods recover the clusters.

\subsection{Two moons}

\begin{figure}

\includegraphics{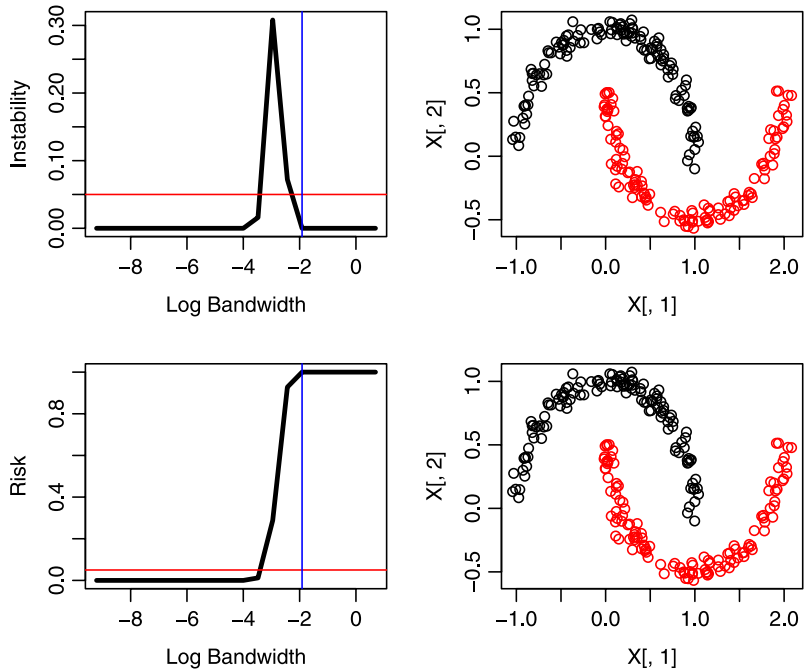}

\caption{Clusters obtained from instability (top)
and excess mass (bottom).
The data are in $\mathbb{R}^{20}$ but only the first two
components are plotted.}
\label{fig::HighDim}
\end{figure}

This is a 20-dimensional example.
The data lie on two half-moons embedded in
$\mathbb{R}^{20}$.
The results are shown in Figure
\ref{fig::HighDim}.
Only the first two coordinates of the
data are plotted.
Again we see that both methods recover the clusters.

%
%
%
%
%
%

\section{Discussion}
\label{section::discussion}

%

As is common in density clustering,
we have assumed a fixed, given value of $\lambda$.
In practice, we recommend that the results should be
computed for a range of values
of $\lambda$ [see, e.g., \citet{stuelzenugent2008}, and
references therein].
It is important to choose a different bandwidth for each $\lambda$.
Indeed, inspection of the proof of Theorem \ref{thm::excessmass}
shows that the optimal bandwidth is a function of $\lambda$
and that
$h(\lambda) \to0$ as $\lambda$ increases.
Further research on data-dependent methods to choose $\lambda$ and
$\rho$
(the parameter used in the graph-based algorithm of Section \ref
{sec::friend}) would be very useful.

We discussed the idea
of using stability to choose a bandwidth.
We saw that the behavior of the selected bandwidth
is quite different than with the excess mass method.
This method seems to work well for density clustering
unlike what happens for $k$-means clustering
[\citet{bendavid2006}].
We believe that the stability method deserves
more scrutiny.
In particular, it would be helpful
to understand the behavior of the stability measure
under more general conditions.
Also, the detailed theoretical properties of the modified method for selecting
$h$ based on stability should be explored.

Finally, we note that
there is growing interest in
spectral clustering methods [\citet{Lux}].
We believe there are connections between the work reported here and
spectral methods.

\section{Proofs}
\label{section::proofs}

\begin{pf*}{Proof of Lemma \ref{lem:ph}}
The weak convergence follows from the fact that $P$ is a Radon
measure [see, e.g., \citet{giovanni07}, Theorem~2.79]. As for the
second part, if $x \in S_i$, where $S_i$ has Hausdorff dimension~$d$,
then $p(x) = \pi_i p_i(x)$, with $p_i$ a Lebesgue density, and the
result follows directly from Leoni and Fonseca (\citeyear{giovanni07}),
Theorem 2.73, part (ii). See also the \hyperref[app]{Appendix}. On the
other hand if
$d_i < d$, then it\vspace*{1pt}
is necessary to modify the arguments as follows. Since $K$ is smooth
and supported\vspace*{-2pt} on $B(0,1)$, there exists a $\eta$ such that $ K (
\frac{x -y}{h} ) > \eta$ if $\| x - y \| <
\eta h$. Set $C = \frac{\eta^{d_i+1} v_{d_i}}{c_d}$, where $v_{d_i}$
is the volume of the unit Euclidean ball in $\mathbb{R}^{d_i}$. Then
\begin{eqnarray*}
p_h(x) & = & \frac{1}{c_d h^d} \int_{S_i \cap B(x,h)} K \biggl( \frac
{x -y}{h} \biggr) \,dP(y)\\
& \geq& \frac{1}{c_d h^d} \eta\int_{S_i \cap B(x,\eta h)} dP(y)\\
& = & \frac{\eta^{d_i+1} v_{d_i}}{c_d h^{d - d_i}} \frac{1}{v_{d_i}
(\eta h)^{d_i} } P_i(B(x,\eta h)) \\
& = & \frac{C}{ h^{(d - d_i)}} \frac{P_i(B(x,\eta h))}{v_{d_i} (\eta
h)^{d_i}}.
\end{eqnarray*}
As as $h \to0$, $\frac{P_i(B(x,\eta,h))}{v_{d_i} (\eta h)^{d_i}} \to
p_i(x) < \infty$, by (\ref{eq:pi}) almost everywhere $\mathcal
{H}^{d_i}$, while $\frac{C}{ h^{(d - d_i)}} \to\infty$, thus showing
that $\lim_{h \to0} p_h(x) = \infty$.
\end{pf*}
\begin{pf*}{Proof of Lemma \ref{lem:clust.numb}}
By assumption (C2), for any $0 \leq\varepsilon< \overline{\varepsilon}$
and $0 < h < \overline{h}$,
\[
N_h(\lambda- \varepsilon) = N_h(\lambda) = N_h(\lambda+ \varepsilon) =
N(\lambda) = k.
\]
On the event $\mathcal{E}_{h,\varepsilon}$ it holds that
\[
L_h(\lambda+ \| p_ h- \widehat{p}_h \|_\infty) \subseteq
\widehat{L}_h(\lambda) \subseteq L_h(\lambda- \| p_ h- \widehat
{p}_h \|_\infty),
\]
which implies that, on the same event,
\[
k = N_h(\lambda+ \| p_h - \widehat{p}_h \|_\infty) \leq\widehat
{N}_h(\lambda) \leq N_h(\lambda- \| p_h- \widehat{p}_h \|_\infty) = k.
\]
\upqed\end{pf*}
\begin{pf*}{Proof of Lemma \ref{lem:xi.theta}}
Recall that $K_h$ is supported on $B(0,h)$. For the first claim, it
is enough to show that, for any $\varepsilon\in[0,
\overline{\varepsilon})$, $L_h(\lambda- \varepsilon) - L(\lambda-
\varepsilon) \subseteq\partial L(\lambda- \varepsilon) +
B(0,h)$. Indeed, by (C3), $ \mu( \partial L(\lambda-
\varepsilon) \oplus B(0,h) ) \leq C_2 h^\theta, $ which
implies~(\ref{eq:theta}). Thus, we will prove that, if $w \notin\partial
L(\lambda- \varepsilon) \oplus B(0,h)$, then $w \notin L_h(\lambda-
\varepsilon) - L(\lambda- \varepsilon)$. For such a point $w$, either
$p(w) \geq\lambda- \varepsilon$ or, by conditions (C2), $p(z) < \lambda
- \varepsilon$ for every $z \in B(w,h)$. Since the kernel $K$ has
compact support, the latter case implies that $p_h(w) < \lambda-
\varepsilon$ as well. Therefore,
\begin{eqnarray*}
w & \in& \{x \dvtx p(x) \geq\lambda- \varepsilon\} \cup\{ x \dvtx
p_h(x) < \lambda- \varepsilon\}\\
& = & \{ x \dvtx p(x) <\lambda- \varepsilon, p_h(x) \geq\lambda-
\varepsilon\}^c \\
& = & \bigl( L_h(\lambda- \varepsilon) - L(\lambda- \varepsilon)
\bigr)^c.
\end{eqnarray*}
As for inequality (\ref{eq:xi}), it is enough to observe that the set
\[
I_{h,\varepsilon} = \bigl( L_h(\lambda- \varepsilon) - L(\lambda-
\varepsilon) \bigr) \cap S
\]
either has zero probability (because it is empty or has Lebesgue
measure $0$) or has positive Lebesgue measure. In the former case, we
obtain $\xi= \infty$. In the latter case, $I_{h,\varepsilon}$ must be
full dimensional, so that, by (\ref{eq:theta}), $\mu(I_{h,\varepsilon})
\leq C_3 h$, for all $h \in(0,\overline{h})$. Since $p$ is bounded by
$\lambda$ on $I_{h,\varepsilon}$, we obtain
\[
P\bigl(L_h(\lambda- \varepsilon) - L(\lambda- \varepsilon)\bigr) =
P(I_{h,\varepsilon
}) \leq\lambda C_2 h =C_3 h,
\]
which implies that we can take $\xi= 1$.
\end{pf*}
\begin{pf*}{Proof of Lemma \ref{lem:1}}
Since $p$ is Lipschitz and integrable, $p^{-1}(\lambda)$ is
$\mathcal{H}^{d-1}$-measurable, so the integral
$\mathcal{H}^{d-1}(\{ x \dvtx p(x) = \lambda\})$ is well defined
for $\lambda\in(0,\| p\|_{\infty})$, where $\mathcal{H}^{d-1}$
denote the $(d-1)$-dimensional Hausdorff measure
in~$\mathbb{R}^d$. Furthermore, we can use the coarea formula. See
\citet{evansgariepy92} and \citet
{ambrosiofuscopallara00} for
backgrounds on Hausdorff measures and the coarea formula. By
the Rademacher theorem, the set $E_1$ of points where $p$ is not
differentiable has Lebesgue measure zero. By Lemma 2.96 in
\citet{ambrosiofuscopallara00}, the set $E_2 = \{ x \dvtx\|
\nabla p(x) \| = 0\}$ is such that $\mathcal{H}^{d-1} \{
p^{-1}(\lambda) \cap E_2 \} = 0$, for all $\lambda\in(0,\|
p\|_{\infty})$ outside of a set $E_3 \subset\mathbb{R}$ of Lebesgue
measure $0$. Without loss of generality, below we may assume that
$E_1$ and $E_2$ are empty. Thus, we can assume that, for any
$\lambda\in(0,\| p\|_{\infty}) \cap E_3^c$, there exists positive
numbers $\overline{\varepsilon}$, $C$ and $M$ such that:

\begin{longlist}
\item$\inf_{x \in\{ x \dvtx|p(x) - \lambda| <
\overline{\varepsilon} \} }\| \nabla p(x) \| > C $, almost
everywhere-$\mu$;
\item$\sup_{ \eta\in(-\overline{\varepsilon}, \overline{\varepsilon}) }
\mathcal{H}^{d-1}(\{ x \dvtx p(x) = \lambda+ \eta\}) < M$.
\end{longlist}
Then for each $\varepsilon\in(0, \overline{\varepsilon})$,
\begin{eqnarray*}
P\bigl(\{ x \dvtx|p(x) - \lambda| < \varepsilon\}\bigr) & = & \int p(x) 1_{\{
|p(x) - \lambda| < \varepsilon\}} \,d \mu(x)\\
& = & \int\frac{p(x)}{\| \nabla p(x) \|} 1_{\{ |p(x) - \lambda| <
\varepsilon\}} \| \nabla p(x) \| \,d \mu(x) \\
& = & \int_{-\varepsilon}^{+\varepsilon} \int_{\{ p^{-1}(\lambda+ u)\}}
\frac{p(x)}{\| \nabla p(x) \|} \,d \mathcal{H}^{n-1}(x)\, du \\
& = & \int_{-\varepsilon}^{+\varepsilon} (\lambda+ u) \int_{\{
p^{-1}(\lambda+ u)\}} ( \| \nabla p(x) \| )^{-1} \,d
\mathcal{H}^{n-1}(x) \,du \\
& \leq& \frac{2 \lambda M}{C} \varepsilon,
\end{eqnarray*}
where the second equality holds because $\| \nabla p(x) \|$ is bounded
away from $0$ on $\{ x \dvtx|p(x) - \lambda| < \varepsilon\}$ by (i),
the third equality is a direct application of the coarea formula
[see, e.g., Proposition 3, page 118 in \citet{evansgariepy92}] and
the last inequality follows from (i) and (ii).
\end{pf*}
\begin{pf*}{Proof of Corollary \ref{cor:2}}
Following the proof of Lemma \ref{lem:1} and using our additional
assumption that $p$ is of class $\mathcal{C}^1$, without any loss of
generality, below we can assume that the set $E_1$ and $E_2$ are
empty and we recall that $E_3$ has Lebesgue measure
$0$.
Let $\lambda\notin E_3$ be such that
\[
{\inf_{x \in p^{-1}(\lambda)}} \| \nabla p(x) \| > 0.
\]
We now claim that there exists a nonempty neighborhood $U$ of $\lambda
$ for which
\[
{\inf_{\lambda\in U} \inf_{x \in p^{-1}(\lambda)}} \| \nabla p(x) \|
> 0.
\]
Indeed, arguing by contradiction, suppose that the previous display
were not verified for any neighborhood $U$ of
$\lambda$.
Then there exist sequences $\{ \lambda_n \} \subset\mathbb{R}$ and
$\{ x_n \} \subset S$ such that $\lim_n \lambda_n = \lambda$, and $x_n
\in p^{-1}(\lambda_n)$ and $\nabla p(x_n) = 0$ for each~$n$. By
compactness, it is possible to extract a subsequence $\{ x_{n_k} \}$
of $\{ x_n \}$ such that $x_{n_k} \rightarrow x$, for some $x \in
p^{-1}(\lambda)$. Since $p$ is of class $\mathcal{C}^1$, this implies
that $\nabla p(x_{n_k}) \rightarrow\nabla p(x)$ as well. However,
$\nabla p(x_{n_k}) = 0$ for each $k$ by construction, while $\nabla
p(x) \neq0$. This produces a contradiction. Thus, for each $\lambda$
that is not a critical point, one can find a neighborhood of positive
length containing it and, by Lemma \ref{lem:1}, (C1) holds at $\lambda$
with $\gamma= 1$. Since, using compactness again, $\| p \|_\infty<
\infty$, this implies that there can only be a finite number of
critical points for which $\gamma$ may differ from $1$.
\end{pf*}
\begin{pf*}{Proof of Lemma \ref{lemma::handy2}}
Since $\varepsilon< \overline{\varepsilon}$ and $h < \overline{h}$, in
virtue of (C2)(b) it holds that, on $\mathcal{E}_{h,\varepsilon}$,
\[
\widehat{L}_h(\lambda) \supseteq L_h(\lambda+ \varepsilon) \supseteq
L(\lambda+ \varepsilon)
\]
and
\[
\widehat{L}_h(\lambda) \subseteq L_h(\lambda- \varepsilon) =
L(\lambda- \varepsilon) \cup\bigl( L_h(\lambda- \varepsilon) -
L(\lambda- \varepsilon) \bigr).
\]
Because $L(\lambda+ \varepsilon) \subseteq L(\lambda) \subseteq
L(\lambda- \varepsilon)$,
the above inclusions imply, still on $\mathcal{E}_{h,\varepsilon}$, that
\begin{eqnarray*}
\widehat{L}_h(\lambda) \Delta L(\lambda) & \subseteq& \bigl(
L(\lambda- \varepsilon) - L(\lambda+ \varepsilon) \bigr) \cup\bigl(
L_h(\lambda- \varepsilon) - L(\lambda- \varepsilon) \bigr)\\
& = & A \cup B,
\end{eqnarray*}
where it is clear that the sets $A$ and $B$ are disjoint. Taking
expectation with respect to $P$ of the indicators of the sets $\widehat
{L}_h(\lambda) \Delta L(\lambda)$, $A$ and $B$ and using condition
(C1) and Lemma \ref{lem:xi.theta} yield (\ref{eq:firstrates}).
\end{pf*}
\begin{pf*}{Proof of Proposition \ref{prop::density}}
The claimed results are a direct consequence of Corollary 2.2 in
\citet{Gine02}.
We outline the details below. We rewrite the left-hand side of (\ref
{eq:sup}) as
\[
\mathbb{P} \Biggl\{ \Biggl\| \sum_{i=1}^n f(X_i) - \mathbb
{E}[f(X_1)] \Biggr\|_{\mathcal{F}_h} > 2 \varepsilon n h^d \Biggr\},
\]
where
\[
\mathcal{F}_h = \biggl\{ K \biggl( \frac{x - \cdot}{h} \biggr), x
\in\mathbb{R}^d \biggr\}
\]
and then proceed to apply Gin\'{e} and Guillou (\citeyear{Gine02}),
Corollary 2.2.
Following their notation, we set $t = n h^d \varepsilon$ and, since,
\[
\sup_{f \in\mathcal{F}_h} \operatorname{Var}[f] \leq\sup_z \int
_{\mathbb{R}^d} K^2 \biggl( \frac{z - x}{h} \biggr) \,dP(x) \leq h^d D,
\]
we can further take $\sigma^2 = h^d D$ and $U = C\| K \|_\infty$,
where $C$ is a positive constant, depending on $h$, such that $\sigma<
U/2$. Then conditions (2.4), (2.5) and (2.6) of \citet{Gine02} are
satisfied for all $n$ bigger than some finite $n_0$, which depends
on the VC characteristics of $K$, $D$, $\| K \|_\infty$, $C$
and $\varepsilon$.
Part 2 is proved in a very similar way. In this case, we take the
supremum over the the entire class $\mathcal{F}$ and we set $\sigma
_n^2 = h_n^d D$ and
$U = \| K \|_{\infty}$. For all $n$ large enough, condition (2.5) is
trivially satisfied because $h_n = o(1)$, while
equations (2.4) and (2.6) hold true by virtue of~(\ref{eq:hn.en}).
The unspecified constants again depend 
on the VC characteristics of $K$, $D$ and $\| K \|_\infty$.
\end{pf*}
\begin{pf*}{Proof of Theorem \ref{thm::levelset}}
We can write
%
%
\begin{equation}\label{eq:level.dec}\qquad
\mathbb{E}(\rho(p,\widehat{p}_h,P)) = \mathbb{E} \biggl( \int
_{\widehat{L}_h(\lambda) \Delta L(\lambda)}dP; \mathcal
{E}_{h,\varepsilon} \biggr) + \mathbb{E} \biggl(\int_{\widehat
{L}_h(\lambda) \Delta L(\lambda)}dP; \mathcal{E}_{h,\varepsilon}^c
\biggr),
\end{equation}
where for a random variable $X$ defined on some probability space
$(\Omega, \mathcal{F},\mathbb{P})$ and an event $\mathcal{E}
\subset\mathcal{F}$, $\mathbb{E}(X;\mathcal{E}) \equiv\int
_{\Omega\cap A} X(\omega) \,d\mathbb{P}(\omega)$.
Using Proposition \ref{prop::density}, the second term on the
right-hand side is upper bounded by
%
%
\begin{equation}\label{eq:E}
\mathbb{P}(\mathcal{E}_{h,\varepsilon}) \leq L e^{- n C_K h^d e^2}.
\end{equation}
As for the first term on the right-hand side of (\ref{eq:level.dec}),
without loss of generality, we consider separately the case in which
the support of $P$ has no lower-dimensional components and the case in
which it of lower dimension. The result for the cases in which the
support has components of different dimensions follows in a
straightforward way.

If the support of $P$ consists of full-dimensional sets, then, on the
event $\mathcal{E}_{h,\varepsilon}$,
\begin{eqnarray*}
\int_{\widehat{L}_h(\lambda) \Delta L(\lambda)}dP & \leq& P
\bigl( L(\lambda- \varepsilon) - L(\lambda+ \varepsilon) \bigr) + P \bigl(
L_h(\lambda- \varepsilon) - L(\lambda- \varepsilon) \bigr)\\
& \leq& C_1 \varepsilon^\gamma+ C_3 h^\xi,
\end{eqnarray*}
where the first inequality stems from (\ref{eq:firstrates}) and the
second from conditions (C1) and (\ref{eq:xi}).

If instead $P$ has lower-dimensional support, then, because, on the
event $\mathcal{E}_{h,\varepsilon}$, $\widehat{L}_h \subset L_h(\lambda
- \varepsilon)$ and because $L \subset L_h(\lambda- \varepsilon)$
by (C2)(b), we see that, on $\mathcal{E}_{h,\varepsilon}$,
\[
\int_{\widehat{L}_h(\lambda) \Delta L(\lambda)}dP = 0.
\]
We conclude that $\mathbb{E}(\rho(p,\widehat{p}_h,P);\mathcal
{E}_{h,\varepsilon})$ is bounded by $\max\{C_1,C_3\} (\varepsilon^\gamma+
h^\xi)$ if the support of $P$ contains a full-dimensional set and is
$0$ otherwise. This, combined with (\ref{eq:E}), yields the claimed
upper bound on the level set risk with $C_L = \max\{C_1,C_3,L \}$.
The convergence rates are established using simple algebra. Notice that
the choice of the sequences $\{ \varepsilon_n\}$ and $\{ h_n\}$ does not
violate condition (\ref{eq:hn.en}).
\end{pf*}
\begin{pf*}{Proof of Corollary \ref{corollary::fI}}
For each $i \in\{1, \ldots,n \}$,
\[
\mathbb{P} ( i \in\widehat{I}_h| \mathcal{E}_{h,\varepsilon
} ) \leq\mathbb{P}
( X_i \in\widehat{L}_h \Delta L | \mathcal{E}_{h,\varepsilon
} ) \leq\max\{C_1,C_3 \}(\varepsilon^\gamma+ h^\xi) \frac
{1}{\mathbb{P}(\mathcal{E}_{h,\varepsilon})},
\]
where the last inequality is due to Lemma
\ref{lemma::handy2}. Thus,
\begin{eqnarray*}
\mathbb{E} (|\widehat{I}_h| ) & \leq&
\sum_{i=1}^n \mathbb{P} ( i \in\widehat{I}| \mathcal
{E}_{h,\varepsilon} ) \mathbb{P}(\mathcal{E}_{h,\varepsilon}) + n
\mathbb{P}(\mathcal{E}^c_h)\\
& \leq& n \bigl( \max\{C_1,C_3 \} (\varepsilon^\gamma+ h^\xi) +
\mathbb{P}(\mathcal{E}^c_h) \bigr)\\
& \leq& C_L ( \varepsilon^\gamma+ h^\xi+ e^{- C_K n h^d \varepsilon
^2} ).
\end{eqnarray*}
\upqed\end{pf*}
\begin{pf*}{Proof of Theorem \ref{thm::excessmass}}
From (\ref{eq:excessmass}), we have
\[
\mathcal{E}(L) - \mathcal{E}(\widehat{L}_h) =
{\int_{ \widehat{L}_h\Delta L}} |p_0 - \lambda| \,d \mu+ P_1(L) -
P_1(\widehat{L}_h),
\]
where $p_0 = \frac{d P_0}{d \mu}$. Since, on the event $\mathcal
{E}_{h,\varepsilon}$, $\widehat{L}h \supset L_h(\lambda+ \varepsilon)$,
we obtain, on the same event,
\[
P_1(L) - P_1(\widehat{L}_h) \leq P_1(L) - P_1(L_h + \varepsilon) = 0,
\]
where the last equality is due to condition (C2)(b). Therefore,
\[
\mathcal{E}(L) - \mathcal{E}(\widehat{L}_h) \leq
{\int_{ \widehat{L}_h\Delta L}} |p_0 - \lambda| \,d \mu.
\]
Just like in the proof of Theorem \ref{thm::levelset}, we treat
separately the case in which the support of $P$ is of lower-dimension
and the case in which it consists of full-dimensional sets.
If the support of $P$ is not of full dimension, then, on $\mathcal
{E}_{h,\varepsilon}$,
\[
\mathcal{E}(L) - \mathcal{E}(\widehat{L}_h) \leq\lambda\mu
(\widehat{L}_h\Delta L) \leq\lambda\mu\bigl(L_h (\lambda- \varepsilon) -
L(\lambda- \varepsilon)\bigr) \leq\lambda C_2 h^\theta
\]
by (\ref{eq:theta}). 
On the other hand, if the support of $P$ has no lower-dimensional
components (so that $p_0 = p$), still on the event $\mathcal
{E}_{h,\varepsilon}$ and using Lemma \ref{lemma::handy2},
%
%
\begin{eqnarray}\label{eq:mass.dec}
{\int_{\{ \widehat{L}_h(\lambda) \Delta L(\lambda) \}}} | p - \lambda
| \,d \mu
&\leq&{\int_{L(\lambda- \varepsilon) -
L(\lambda+ \varepsilon) }} | p - \lambda| \,d \mu\nonumber\\[-8pt]\\[-8pt]
&&{} + {\int_{L_h(\lambda-
\varepsilon) - L(\lambda- \varepsilon)} }| p - \lambda| \,d
\mu.\nonumber
\end{eqnarray}

The first term on the right-hand side of the previous inequality can be
bounded as follows:
\begin{eqnarray*}
{\int_{L(\lambda- \varepsilon) - L(\lambda+ \varepsilon) } }| p - \lambda
| \,d \mu(x) & = &
{\int_{\{ x \dvtx|p(x) - \lambda| < \varepsilon\} }} | p - \lambda| \,d
\mu(x)\\
& \leq& \varepsilon\int_{\{ x \dvtx|p(x) - \lambda| < \varepsilon\} }
d \mu(x) \\
& = & \frac{\varepsilon}{\lambda- \varepsilon}\int_{\{ x \dvtx|p(x) -
\lambda| < \varepsilon\} }( \lambda- \varepsilon) \,d \mu\\
& \leq& \frac{\varepsilon}{\lambda- \varepsilon} \int_{\{ x \dvtx
|p(x) - \lambda| < \varepsilon\} } p(x) \,d \mu(x)\\
& \leq& \frac{C_1}{\lambda- \varepsilon} \varepsilon^{\gamma+1},
\end{eqnarray*}
where the last inequality is due to condition (C1). As for the
second term of the right-hand side of (\ref{eq:mass.dec}),
\[
{\int_{L_h(\lambda- \varepsilon) - L(\lambda- \varepsilon)} }| p - \lambda
| \,d \mu\leq\lambda\mu\bigl( L_h(\lambda- \varepsilon) - L(\lambda-
\varepsilon) \bigr) \leq\lambda C_2 h^\theta
\]
by (\ref{eq:theta}).


Thus, we conclude that $\mathbb{E} ( \mathcal{E}(L) - \mathcal
{E}(\widehat{L}_h) ; \mathcal{E}^{h,\varepsilon} )$ is bounded
by $\lambda C_2 h^\theta$ if the support of $P$ is a lower-dimensional
set and by
\[
\max\biggl\{ \lambda C_2 , \frac{C_1}{\lambda- \varepsilon} \biggr\}
( \varepsilon^{\gamma+1} + h^\theta)
\]
otherwise. Next, by compactness of $S$, and using (\ref{eq:E}),
\[
\mathbb{E} \bigl( \mathcal{E}(L) - \mathcal{E}(\widehat{L}_h) ;
\mathcal{E}^c_{h,\varepsilon} \bigr) \leq\bigl( 1 + \lambda\mu\bigl(S +
B(0,\overline{h})\bigr) \bigr) \mathbb{P}(\mathcal{E}_{h,\varepsilon
})\leq C_S (1 + \lambda) L e^{- n C_K h^d e^2}
\]
for some positive constant $C_S$, uniformly in $h < \overline{h}$.
The claimed upper bound on the excess mass risk now follows by taking
$C_M = \max\{ \lambda C_2 , \frac{C_1}{\lambda- \varepsilon},C_S
(1 + \lambda) L \}$. The convergence rates can be easily
obtained by simple algebra. Notice that the choice of the sequences $\{
\varepsilon_n\}$ and $\{ h_n\}$ does not violate condition (\ref{eq:hn.en}).
\end{pf*}
\begin{pf*}{Proof of Theorem \ref{thm::split}}
This follows by combining the version of Talagrand's inequality for
empirical processes as given in \citet{massart2000} with an
adaptation of the arguments used in the proof of Theorem 7.1 in
\citet{gyorfi2002}. For completeness, we provide the details.

Define
$\widehat{h} = \operatorname{arg}\operatorname{sup}_{h \in\mathcal
{H}}\widehat
{\mathcal{E}}(L_h)$, where
\[
\widehat{\mathcal{E}}(L_h) = \frac{1}{n} \sum_{i=1}^n I(Z_i \in
L_h) - \lambda\mu(L_h)
\]
and $h_* = \operatorname{arg}\operatorname{sup}_{h \in\mathcal
{H}}\mathcal{E}_X(L_h)$.
Set $\Gamma(h) = \mathcal{E}_X(L_{h_*}) - \mathcal{E}_X(L_{h})$,
where $h \in\mathcal{H}$.
Recall that both $L_{h^*}$ and $L_h = \{ x \dvtx\widehat{p}_h \geq
\lambda\}$,
are random sets depending on the training set~$X$.
We will bound $\mathbb{E}( \Gamma(\widehat{h}) )$, where the
expectation is over the joint distribution of $X$ and $Y$.

We can write
\[
\mathbb{E}(\Gamma(\widehat{h}) |X) =
\underbrace{\mathbb{E}(\Gamma(\widehat{h}) |X) - (1+\delta
)\widehat{\Gamma}(\widehat{h})}_{T_1}{} + {}
\underbrace{(1+\delta)\widehat{\Gamma}(\widehat{h})}_{T_2},
\]
where
$\widehat{\Gamma}(h) = \widehat{\mathcal{E}}(L_{h_*}) - \widehat
{\mathcal{E}}(L_h)$. 
Note that
\[
\widehat\Gamma(\widehat{h}) = \widehat{\mathcal{E}}(L_{\widehat
{h}}) - \widehat{\mathcal{E}}(L_{h_*}) \leq\widehat{\mathcal
{E}}(L_{h_*}) - \widehat{\mathcal{E}}(L_{h_*}) =0.
\]
Thus,
$\mathbb{E}(T_2|X) \leq0$.
We conclude that
%
%
\begin{equation}\label{eq::iter}\qquad
\mathbb{E}(\Gamma(\widehat{h})) =
\mathbb{E}(\mathbb{E}(\Gamma(\widehat{h})|X)) =
\mathbb{E} ( \mathbb{E}(T_1|X)) + \mathbb{E} ( \mathbb{E}(T_2|X))
\leq
\mathbb{E} ( \mathbb{E}(T_1|X)).
\end{equation}

Now we bound
$\mathbb{E}(T_1|X)$.
Consider the empirical process
\[
Z = \sup_{h \in\mathcal{H}} \widehat{\Gamma}(h),
\]
so that $Z = \widehat{\Gamma}(\widehat{h})$ and $\mathbb{E}(\Gamma
(\widehat{h})|X) = \mathbb{E}(Z|X)$. We have
\begin{eqnarray*}
\mathbb{P}(T_1 \geq s |X) &=&
\mathbb{P} \bigl(
\mathbb{E}(Z|X) - (1+\delta) Z \geq
s \mid X \bigr)\\
& = &
\mathbb{P} \biggl( \mathbb{E}(Z|X) - Z \geq
\frac{s + \delta\mathbb{E}(Z|X)}{1+\delta} \Bigm| X \biggr).
\end{eqnarray*}
Notice that, conditionally on $X$, $Z = \frac{1}{n} \sup_{h \in
\mathcal{H}} \sum_{i=1}^n f_h(Y_i)$, where, for each \mbox{$h \in\mathcal
{H}$}, $f_h \dvtx\mathbb{R}^d \mapsto\mathbb{R}$ is the function
given by
\[
f_h(x) = I(x \in L_{h_*}) - \lambda\mu(L_{h_*}) - \bigl(I(x \in L_h)
- \lambda\mu(L_h) \bigr)
\]
with $\| f_h \|_{\infty} < \kappa$.
Let $\sigma^2 \equiv\mathbb{E}( \frac{1}{n} \sup_{h \in\mathcal{H}}
\sum_{i=1}^n f^2_h(Y_i) |X)$ and notice that $\sigma^2 \leq\kappa
\mathbb{E}(\sup_h \widehat{\Gamma}(h) |X) = \kappa\mathbb{E}(Z|X)$.
Thus,
\[
\mathbb{P}(T_1 \geq s |X) \leq\mathbb{P} \biggl( \mathbb{E}(Z|X) - Z
\geq
\frac{s + \delta\sigma^2 / \kappa}{1+\delta} \Bigm| X \biggr),
\]
which, by Corollary 13 in \citet{massart2000}, is upper bounded by
\[
2 \exp\biggl\{ - \frac{n ( ({s + \delta\sigma^2 / \kappa
})/({1+\delta}) )^2}{4(4 \gamma^2 \sigma^2 + {7}/{4} \kappa
\varepsilon) } \biggr\}.
\]
Then, some algebra [see Problem 7.1 in \citet{gyorfi2002}]
yields the final bound
\[
\mathbb{P}(T_1 \geq s |X) \leq2 \exp\biggl\{ \frac{-n s}{d(\delta
,\kappa)} \biggr\},
\]
where $d(\delta,\kappa)$ is given the in the statement of the theorem.

Set $u = \frac{d(\delta,\kappa)}{n}\log2$. Then
\begin{eqnarray*}
\mathbb{E}(T_1|X) &=&
\int_0^\infty\mathbb{P}(T_1>s|X) \,ds \leq
u + \int_u^\infty\mathbb{P}(T_1>s|X) \,ds \\
&=&
u + \frac{2 d(\delta,\kappa) }{n} \exp\biggl\{ - \frac
{nu}{d(\delta,\kappa)} \biggr\}\\
& = &
d(\delta,\kappa) \frac{1+\log2}{n}.
\end{eqnarray*}
From (\ref{eq::iter}), we conclude that
\[
\mathbb{E}(\Gamma(\widehat{h})) \leq d(\delta,\kappa) \frac
{1+\log2}{n}
\]
and so
\[
\mathbb{E}(M(\widehat{h})) \leq\mathbb{E}(M(h_*)) + d(\delta
,\kappa) \frac{1+\log2}{n},
\]
which implies that
\[
\mathbb{E}(\mathcal{E}(\widehat{h})) \geq\mathbb{E}(\mathcal
{E}(h_*)) - d(\delta,\kappa) \frac{1+\log2}{n}.
\]
This shows (\ref{eq::expect}).
\begin{pf*}{Proof of Theorem \ref{thm::adapt}}
Define
$r_n(\gamma,\theta) = (\frac{\log n}{n} )^{{\theta
(\gamma+1)}/({2\theta+ d(\gamma+1)})}$.
For each $\theta$,
$r_n(\gamma,\theta)$ is decreasing in $\theta$ and
\[
r_n(\Upsilon_n(\theta),\theta) \leq2 r_n(\infty,\theta).
\]
Hence,
$\inf_{\gamma\in[0,\Upsilon_n(\theta)]} r_n(\gamma,\theta) \leq
2\inf_{\gamma\geq0} r(\gamma,\theta)$.
Some algebra shows that
$|\partial r_n(\gamma,\break\theta)/\partial\gamma| \leq A_n(\theta)$
for all $\gamma$ and $\theta$.
Therefore, for each $j$,
$r_n(\gamma_j(\theta),\theta) = r_n(j\delta_n(\theta) + \delta
_n(\theta),\theta)\geq
r_n(j \delta_n(\theta),\theta) - \delta_n(\theta)A_n(\theta) \geq
r_n(\gamma_j(\theta),\theta)/2$.
Let $h_n = h(\gamma,\theta)$.
By Theorem \ref{thm::excessmass},
$R^M(p,\widehat{p}_{h_n}) = O((\log n /n)^{{\theta(\gamma
+1)}/({2\theta+ d(\gamma+1)})}$.
Let $h_*\in\mathcal{H}_n$ minimize
$R^M(p,\widehat{p}_{h})$ for $h\in\mathcal{H}_n$.
Then
$R^M(p,\widehat{p}_{h_*}) \leq2 R^M(p,\widehat{p}_{h_n})$.
So,
\begin{eqnarray*}
R^M(p,\widehat{p}_{\widehat h}) & \leq& d(\delta,\kappa)\frac{1 +
\log2}{n} + R^M(p,\widehat{p}_{h_*})\\
& \leq& d(\delta,\kappa)\frac{1 + \log2}{n} + 2 R^M(p,\widehat
{p}_{h_n})\\
&=& d(\delta,\kappa)\frac{1 + \log2}{n} + 2 r_n(\gamma,\theta)\\
&=&O \biggl(\frac{\log n}{n} \biggr)^{{\theta(\gamma
+1)}/({2\theta+ d(\gamma+1)})}.
\end{eqnarray*}
\upqed\end{pf*}
\noqed\end{pf*}
\begin{pf*}{Proof of Theorem \ref{thm::stability}}
(1)
When $h=0$,
$\{\widehat{p}_h >\lambda\} = X$ and
$\{\widehat{q}_h >\lambda\} = Y$
so that
$\{\widehat{p}_h >\lambda\}\Delta\{\widehat{q}_h >\lambda\} =(X,Y)$.
Since $P$ has a Lebesgue density, with probability one,
$d\widehat{P}_Z$ puts no mass on $(X,Y)$ and, therefore,
$\Xi(0) =0$. By compactness of $S$,
if $h \geq\operatorname{diam}(S)$, then
$ \| \widehat{p}_h\|_\infty= \|\widehat{q}_h\|_\infty= \frac
{1}{h^d v_d}$,
with the supremum attained by any $z \in S$. Thus, as $h\to\infty$,
$\|\widehat{p}_h - \widehat{q}_h \|_\infty\to0$ and consequently,
$\Xi(\infty) \to0$.

(2) Note that
\begin{eqnarray*}
\Xi(h) &=& \rho(\widehat{p}_h,\widehat{q}_h,\widehat{P}_Z) =
\int_{ \{\widehat{p}_h \geq\lambda\}\Delta\{\widehat{q}_h \geq
\lambda\}}d \widehat{P}_Z(z)\\
&=&
\int I\bigl(\widehat{p}_h(z) \geq\lambda,\widehat{q}_h(z) \leq\lambda\bigr)
\,d\widehat{P}_Z(z) \\
&&{} +
\int I\bigl(\widehat{p}_h(z) \leq\lambda,\widehat{q}_h(z) \geq\lambda\bigr)
\,d\widehat{P}_Z(z).
\end{eqnarray*}
Define $\xi(h) = \mathbb{E}(\Xi(h)|X,Y)$. Then
\begin{eqnarray*}
\xi(h) &=& \rho(\widehat{p}_h,\widehat{q}_h,P)\\
&=&
\int I\bigl(\widehat{p}_h(z) \geq\lambda,\widehat{q}_h(z) \leq\lambda\bigr)
\,dP(z) \\
&&{}+
\int I\bigl(\widehat{p}_h(z) \leq\lambda,\widehat{q}_h(z) \geq\lambda\bigr)
\,dP(z)\\
&\stackrel{d}{=}&
2 \int I\bigl(\widehat{p}_h(z) \geq\lambda,\widehat{q}_h(z) \leq\lambda
\bigr) \,dP(z),
\end{eqnarray*}
where $\stackrel{d}{=}$ denotes identity in distribution. Let
$\pi_h(z) = \mathbb{P}(\widehat{p}_h(z) \leq\lambda) = \mathbb
{P}(\widehat{q}_h(z) \leq\lambda) $.
By Fubini's theorem and independence,
%
%
\begin{eqnarray}\label{eq::Exi}
\mathbb{E}(\Xi(h)) &=& \mathbb{E}( \xi(h)) \nonumber\\
& = & 2 \int_{\mathbb{R}^d} \mathbb{P}\bigl(\widehat{p}_h(z) \geq
\lambda,\widehat{q}_h(z) \leq\lambda\bigr) \,dP(z)
\nonumber\\[-8pt]\\[-8pt]
&=&
2 \int_{\mathbb{R}^d} \mathbb{P}\bigl(\widehat{p}_h(z) \geq\lambda
\bigr)\mathbb{P}\bigl(\widehat{q}_h(z) \leq\lambda\bigr) \,dP(z) \nonumber\\
&=&
2 \int_{\mathbb{R}^d} \pi_h(z)\bigl(1-\pi_h(z)\bigr) \,dP(z).\nonumber
\end{eqnarray}
Since $\pi_h(z)(1 - \pi_h(z)) \leq1/4$ for all $n$, $h$ and $z$,
(2) follows.

(3)
Let $W = (X,Y)$ be the $2n$-dimensional vector obtained by
concatenating $X$ and $Y$ and define the event
\[
\mathcal{A}_h = \{ B(W_i, h) \cap B(W_j, h) = \varnothing,
\forall i \neq j \}.
\]
Let $h$ be small enough such that $\lambda n h^d v_d < 1$
(trivially satisfied if $\lambda= 0$). Then, for any realization $w$
of the vector $W$ for which the event $\mathcal{A}_h$ occurs,
\[
\int I\bigl(\widehat{p}_h(z) \geq\lambda, \widehat{q}_h(z) \leq\lambda
\bigr)\,d P(z) = \sum_{i=1}^{2n} P(B(w_i,h)).
\]
By our assumptions,
\[
2n \delta h^d v_d \leq\sum_{i=1}^{2n} P(B(w_i,h)) \leq2n \overline
{\Delta} h^d v_d.
\]

Using the union bound, we also have
\[
\mathbb{P} ( \mathcal{A}_h^c ) \leq\pmatrix{2n \cr2}
(2h)^d v_d \overline{\Delta}.
\]
Thus, it follows that, for fixed $n$, $\mathbb{E} ( \xi(h)
) \rightarrow0$ as $h \rightarrow0$ according to
\[
2n \delta h^d v_d \leq\mathbb{E}
( \xi(h) ) \leq h^d v_d 2 \overline{\Delta} \max
\{ 2^d n(2n-1), 2n \}.
\]

(4)
By the same arguments used in the proof of point (1), for all $h \geq h_*$,
$\xi(h) = 0$ almost everywhere with respect to the joint distribution
of $X$ and $Y$, and, therefore,
$\mathbb{E}(\xi(h)) = 0$.
Thus, we need only to consider the case $0 < h \leq h_*$.

Set $p_{z,h} = P(B(z,h))$ and denote with $X_{z,h}$ a
random variable with distribution $\operatorname{Bin}(n, p_{z,h})$. Then
\[
\mathbb{P}\bigl(\widehat{p}_h(z) = 0\bigr) = \mathbb{P}(X_{z,h} = 0) = (
1 - p_{z,h} )^n.
\]
For each $z \in S$, set $D(z,h) = \{ z' \in S \dvtx\| z - z'\| <
h\}$ and $S_h = \{ z \dvtx D(z,h) \neq S\}$. Furthermore, set
$p_{h,\max} = \sup_{z \in S_h} \{ p_{z,h}\}$ and $p_{h,\min} = \inf_{z
\in S_h} \{ p_{z,h}\}$. Then the expected instability can be
written as
\[
\mathbb{E} ( \Xi(h) ) = 2 \int_{S_h} \pi_h(z) \bigl(1 - \pi
_h(z)\bigr) \,d P(z)
\]
so that $A_{h,n} \leq\mathbb{E} ( \Xi(h) ) \leq
B_{h,n}$, where
\begin{eqnarray*}
A_{h,n} & \equiv& 2 P(S_h) (1 - p_{h,\max})^n \bigl( 1- (1 -
p_{h,\min})^n \bigr),\\
B_{h,n} & \equiv& 2 P(S_h) (1 - p_{h,\min})^n \bigl( 1- (1 -
p_{h,\max})^n \bigr).
\end{eqnarray*}

We will now upper bound $B_{h,n}/2 $. For the first term, we proceed as follows.
There exists a sphere
$E=B(z_0,h_*/2)$
such that
$ S\subset E$.
[E.g., choose any two points $z,z'$
such that $\Vert z-z'\Vert=h^*$.
Set $z_0 = (z+z')/2$.]
Let
$A = B(z_0,h_*/2) - B(z_0,(h_*-h)/2)$.
We claim that
$S_h\subset A$.
This follows since if
$z\in A^c\cap S$ then $z\in B(z_0,h/2)$ and then
$\sup_{z'\in S}\Vert z-z'\Vert\leq\sup_{z\in B(z_0,h/2),z'\in
B(z_0,h_*/2)}\Vert z-z'\Vert = h$.
Thus, if $z\in S_h$ then
$z\in A\cap S \subset A$.
Hence,
\[
P(S_h) \leq P(A) \leq\overline{\Delta}\mu(A) =
\overline{\Delta}
\frac{ ((h_*/2)^d - (h/2)^d)\pi^{d/2}}{\Gamma((d/2) + 1)} \leq D_1
(h_*-h),
\]
where
\[
D_1 = \frac{\pi^{d/2} h_*^{d-1}}{2^d \Gamma((d/2)+1)}.
\]

For the second term,
let $z_0 = \arg\min_z p_{z,h}$.
Then
\begin{eqnarray*}
1 - p_{h,\min} &=& 1 - P(B(z_0,h))= P(B(z,h_*)) - P(B(z_0,h))\\
&=& P\bigl(B(z,h_*)- B(z_0,h)\bigr) \leq\overline{\Delta} \mu\bigl(B(z,h_*)-
B(z_0,h)\bigr)\\
& \leq& \overline{\Delta}
\frac{ (h_*^d - h^d)\pi^{d/2}}{\Gamma((d/2) + 1)} = D_2 (h_*-h),
\end{eqnarray*}
where
$D_2 = \frac{\pi^{d/2} h_*^{d-1}}{\Gamma((d/2)+1)}$.
The third term is bounded above by 1.
Hence,
$B_n \leq D_1 D_2^n (h_*-h)^{n+1}$.

Now we lower bound $A_{h,n}/2$.
First, we claim that
$S_h$ contains the intersection of a sphere of radius $r/2$
where $r = h_*-h$, with $S$.
Indeed, let $z\in S_h$.
Then there exists
$z'\in S$ such that
$\Vert z-z'\Vert \leq h_* = h+r$.
Let $w\in B(z',r/2)$.
By the triangle inequality,
$\Vert w-z\Vert \geq h+ r/2$.
So $B(z',r/2)\cap S \subset S_h$.
Therefore,
\begin{eqnarray*}
P(S_h) &\geq& P\bigl(B(z',r/2)\cap S\bigr) \geq\underline{\Delta}\mu
\bigl(B(z',r/2)\cap S\bigr)\\
& \geq& \delta\underline{\Delta}\mu\bigl(B(z',r/2)\bigr) = D_3
(h_*-h)^d,
\end{eqnarray*}
where
$D_3 = \frac{\delta\underline{\Delta}\pi^{d/2} }{\Gamma((d/2)+1)}$.

To lower bound the second term,
Let $z_0 = \arg\max_z p_{z,h}$.
Then
\begin{eqnarray*}
1 - p_{h,\max} &=& 1 - P(B(z_0,h))= P(B(z,h_*)) - P(B(z_0,h))\\
&=& P\bigl(B(z,h_*)- B(z_0,h)\bigr) \geq\underline{\Delta} \mu\bigl(\bigl(B(z,h_*)-
B(z_0,h)\bigr)\cap S\bigr)\\
& \geq& \underline{\Delta} \delta\mu\bigl(B(z,h_*)- B(z_0,h)\bigr)
= \underline{\Delta} \delta\frac{ (h_*-h)^d \pi^{d/2}}{\Gamma
(d/2+1)}\\
&=& D_4 (h_* - h)^d,
\end{eqnarray*}
where
$D_4 = \frac{\underline{\Delta} \delta\pi^{d/2}}{\Gamma(d/2+1)}$.
Thus,
$(1 - p_{h\max})^n \geq D_4^n (h_* - h)^{nd}$.
For the third term, argue as above that
$1 - p_{h,\min} \leq D_2 (h_*-h)$
so the third term is larger than $1/2$ when
$h$ is close enough to $h_*$.
Hence,
$A_n \geq\frac{D_3}{2} D_4^n (h_*-h)^{d(n+1)}$.
\end{pf*}
\begin{pf*}{Proof of Theorem \ref{thm::friends}}
By our assumptions (see Section \ref{sec:settings}),
\[
0 < \lim_{r \rightarrow0} \frac{P(B(x,r))}{r^{d_i}} < \infty,
\]
where $d_i = \operatorname{dim}(S_i)$, for any $x$ outside of a set of $P_i$
measure zero. By Theorem 5.7 in \citet{mattila99}, $d_i$ is also the
box-counting dimension of $S_i$. Thus, $d^* = \max_i d_i$. Combined
with (\ref{eq:standard}) this implies that, without loss of
generality, we can assume that there exist constants $\overline{C}> 0$
and $\overline{\rho} > 0$ such that for every ball $B$ of radius
$\rho
< \overline{\rho}$ and center in $L(\lambda)$, $P(B) > \overline{C}
\rho^{d^*}$.

Let $\mathcal{A}$ be a covering of $L(\lambda)$ with balls of radius
$\rho/2$ and centers in $L(\lambda)$, with $\rho< \overline{\rho}$.
By compactness of $L$, $|\mathcal{A}| \leq\overline{M} \rho^{-d^*}$,
where $\overline{M}$ depends on $d^*$ and $L(\lambda)$ but not on
$\rho$.

Next, by Lemma \ref{lem:clust.numb}, on the event
$\mathcal{E}_{h,\varepsilon} = \{ \| p_h - \widehat{p}_h \|_\infty<
\varepsilon\}$, the set $\widehat{L}_h$ consists of $k$ disjoint connected
sets. Since $\rho< \delta/2$, this implies, on the same event, that
$\widehat{N}^{G}_h(\lambda) \geq k$. Thus, on the event
$\mathcal{E}_{h,\varepsilon}$, for some $\varepsilon< \varepsilon_1$ to be
specified below, a~sufficient condition for the event
$\mathcal{O}_{h,n}$ to be verified is that every $A \in\mathcal{A}$
contains at least one point from the set $\widehat{J}_{h} \equiv\{ i
\dvtx\widehat{p}_{h}(X_i) \geq\lambda\}$ [similar arguments are used
also in \citet{cuevas2000}, \citet{cadre2007}]. We conclude that the
probability of $\mathcal{O}^c_{h,n}$ is bounded from above by
\[
\mathbb{P}( \mathcal{E}_{h,\varepsilon}^c ) +
\overline{M} \rho^{-d^*} \sup_{A \in\mathcal{A}_{h}}
\mathbb{P} ( \{ X_i \notin A, \forall i \in\widehat
{J}_{h} \}
\cap\mathcal{E}_{h,\varepsilon} ).
\]
Since, on the event $\mathcal{E}_{h,\varepsilon}$ the set
$J_{h_n} = \{ i \dvtx p_{h_n}(X_i) \geq\lambda+ \varepsilon\}$ is
contained in $\widehat{J}_{h} $,
we further have that, for each $A \in\mathcal{A}_{h}$,
%
%
\begin{equation}\label{eq:ineq}
\mathbb{P} ( \{ X_i \notin A, \forall i \in
\widehat{J}_{h} \} \cap\mathcal{E}_{h,\varepsilon} )
\leq\bigl( 1 - P (A \cap\{ p_{h} \geq\lambda+ \varepsilon\}) \bigr)^{n},
\end{equation}
where the inequality stems from the identity among events
\[
\{ X_i \notin A, \forall i \in J_{h} \} =
\bigcap_{i} \bigl\{\bigl\{ \{ p_{h_n}(X_i) \geq\lambda+ \varepsilon\}\cap
A^c \bigr\} \cup\{ p_{h_n}(X_i) < \lambda+ \varepsilon\} \bigr\},
\]
and the independence of the $X_i$'s. By Lemma \ref{lem:suff.cond}, for
any fixed $0 < \tau< 1/2$, there exists a point $y \in L(\lambda)
\cap L_h(\lambda+ \varepsilon)$ such that $B ( y,\frac{\tau\rho
}{2} ) \subset A \cap L_h( \lambda+ \varepsilon)$, for all
$\varepsilon< \varepsilon(\rho,\tau)$. Thus,
\[
P \bigl(A \cap L_h( \lambda+ \varepsilon) \bigr) \geq P \biggl( B \biggl( y,\frac
{\tau\rho}{2} \biggr) \biggr) \geq\overline{C} \biggl(\frac{\tau
\rho}{2} \biggr) ^{d^*}
\]
for all $\varepsilon< \varepsilon(\rho, \tau)$,
where the second inequality is verified since $\frac{\rho\tau}{2} <
\overline{\rho}$.
Set $\varepsilon(\rho) = \min\{ \varepsilon_1, \varepsilon(\rho, \tau)\}
$. The result now follows from collecting all the terms and the
inequality $(1-x)^{n} \leq e^{-n x}$,
valid for all $0 \leq x \leq1$.
\end{pf*}
\begin{pf*}{Proof of Theorem \ref{thm::bootstrap}}
Let $\mathcal{A}_h$ be a covering of $L_h(\lambda)$ by balls of
radius $\rho/2$ and centers in $L_h(\lambda)$. By the same arguments
used in the proof of the Theorem~\ref{thm::friends},
the probability of the event $(\mathcal{O}^*_{h,n})^c$ is bounded by
\[
\mathbb{P}( \mathcal{E}_{h,\varepsilon}^c ) +
\overline{M} \rho^{-d} \sup_{A \in\mathcal{A}_{h}}
\mathbb{P} ( \{ X^*_j \notin A, \forall j \}
\cap\mathcal{E}_{h,\varepsilon} ),
\]
where the probability is over the original sample
$X = (X_1,\ldots,X_n)$ and the bootstrap sample $X^* = (X^*_1,\ldots
,X_N^*)$. Here, the value of $\varepsilon< \varepsilon_1$ used in the
definition of the event $\mathcal{E}_{h,\varepsilon}$ is to be specified
below. Because of compactness of the support of $P$, $\overline{M}$ is
a constant depending only $d$ and $S + B(0,\overline{h})$.

For a set $S \subseteq\mathbb{R}^d$, we denote with $S^{\otimes n}$
the $n$-fold Cartesian product of $S$ and with $P^h_{X^*|X=x}$ the
conditional distribution of the bootstrap sample $X^*$ given $X = x$,
with $x = (x_1,\ldots,x_n)$.
Let $\mathcal{E}_n = \{ x \in
S^{\otimes n} \dvtx\| p_{h}- \widehat{p}_{h} \|_\infty\leq
\varepsilon\}$, where $\widehat{p}_h$ is the kernel density estimate
based on $x$.
Then, for each $A \in\mathcal{A}_{h}$,
\[
\mathbb{P} ( \{ X^*_j \notin A, \forall j \}
\cap\mathcal{E}_{h,\varepsilon} ) =
E_X ( P_{X^*|X} ( (A^c)^{\otimes n} ); \mathcal{E}_n
),
\]
where, if $X \sim P$, $E_X(f(X);\mathcal{E}) \equiv\int_{\{ x \in
\mathcal{E} \}} f(x) \,dP(x)$.
For every $x \in\mathcal{E}_n$, by the conditional independence of
$X^*$ given $X = x$,
\begin{eqnarray*}
P_{X^*|X=x} ( (A^c)^{\otimes n} )& = &
\biggl( 1 - \frac{\int_{A \cap L_h(\lambda)} \widehat{p}_h(v)
\,dv}{\int_{\{ \widehat{p}_h \geq\lambda\}} \widehat{p}_{h}(v)
\,dv} \biggr)^N \\
& \leq& \biggl( 1 - \frac{\int_{A \cap L_h(\lambda+ \varepsilon)}
(p_h - \varepsilon) \,d\mu}{V(h,\varepsilon)} \biggr)^N,
\end{eqnarray*}
where
\[
V(h,\varepsilon) = \int_{L_h(\max\{\lambda- \varepsilon,0\})} (p_h +
\varepsilon) \,d \mu.
\]
By Lemma \ref{lem:suff.cond}, for any fixed $\tau< 1/2$ and each $h$,
there exists a point $y \in L_h(\lambda) \cap L_h(\lambda+ \varepsilon
)$ such that $B ( y,\frac{\tau\rho}{2} ) \subset A \cap
L_h( \lambda+ \varepsilon)$, for all $\varepsilon< \varepsilon(\rho,\tau
)$. Thus,
\begin{eqnarray*}
\int_{A \cap L_h(\lambda+ \varepsilon)} (p_h - \varepsilon) \,d\mu
&\geq&
\int_{B ( y,{\tau\rho}/{2} )} (p_h - \varepsilon)
\,d\mu\\
&=& \int_{B ( y,{\tau\rho}/{2} )} p_h \,d \mu-
\varepsilon\mu\biggl( B \biggl( y,\frac{\tau\rho}{2} \biggr) \biggr).
\end{eqnarray*}
Next,
\[
V(h,\varepsilon) = \int_{L_h(\lambda)} p_h \,d \mu+ \varepsilon\mu\bigl(
L_h(\max\{ \lambda- \varepsilon, 0 \}) \bigr) + \int_{L_h(\lambda)
- L_h(\max\{ \lambda- \varepsilon, 0 \})} p_h \,d \mu.
\]
Following the proof of Lemma \ref{lem:suff.cond}, one can verify that,
because of assumption (G),
$
\inf_{h \in(0, \overline{h})} \mu( L_h(\lambda) - L_h(\max
\{ \lambda- \varepsilon, 0 \}) ) \rightarrow0
$,
as $\varepsilon\rightarrow0$. Thus,
\[
\frac{\int_{A \cap L_h(\lambda+ \varepsilon)} (p_h - \varepsilon) \,d\mu
}{V(h,\varepsilon)} \geq\frac{\int_{B ( y,{\tau\rho}/{2}
)} p_h \,d \mu}{\int_{L_h(\lambda)} p_h \,d \mu} \bigl(1 + o(1)\bigr)
\]
as $\varepsilon\rightarrow0$. Then using (\ref{eq:cond.boot}) and the
facts $\tau< 1/2$ and $\int_{L_h(\lambda)} p_h \,d \mu\leq1$ for
each $h$, we conclude that there exists a $\varepsilon(\rho,\tau)$ such that
\[
\frac{\int_{A \cap L_h(\lambda+ \varepsilon)} (p_h - \varepsilon) \,d\mu
}{V(h,\varepsilon)} \geq C \rho^d
\]
for all $0 < \varepsilon< \varepsilon(\rho,\tau)$ and for some
appropriate constant $C$, independent of $\rho$ and~$h$. Thus,
\[
P_{X^*|X=x} ( (A^c)^{\otimes n} ) \leq e^{ - N C \rho^d }
\]
and the results now follows by setting $\varepsilon(\rho) = \min\{
\varepsilon_1, \varepsilon(\rho,\tau) \}$.
\end{pf*}
\begin{lemma}\label{lem:suff.cond}
Assume conditions \textup{(C2)} and \textup{(G)}. Then, for any $0 < \tau< 1$
and $\rho> 0$, there exists a positive number $\varepsilon(\rho,\tau)$
such that, for all $\varepsilon< \varepsilon(\rho,\tau)$,
%
%
\begin{equation}\label{eq:final}
\sup_{h \in(0,\overline{h})} \sup_{x \in L(\lambda)}
\operatorname{dist}\bigl(x,L_h(\lambda+ \varepsilon)\bigr) < \tau\rho.
\end{equation}
\end{lemma}
\begin{pf}
The claim follows by minor modifications of the arguments used in the
Appendix of \citet{cadre2007}. We provide some details for
completeness and refer to \citet{lee03} for background. Because of
assumption (G) and in virtue of the regular level set theorem
[see, e.g., \citet{lee03}, Corollary~8.10], for any $\varepsilon
\in(0,
\varepsilon_1)$ and $h \in(0,\overline{h})$, the set $\{ x \dvtx p_h(x)
= \lambda+ \varepsilon\}$ is a closed embedded submanifold of
$\mathbb{R}^d$. Let $r(\varepsilon,h)$ be the maximal radius of the
tubular neighborhood around $\{ x \dvtx p_h(x) = \lambda+ \varepsilon
\}$. Set $\overline{r}_h = \inf_{\varepsilon< \varepsilon_1}
r(\varepsilon,h)$
and notice that $\overline{r}_h > 0$ is positive for each $h \in
(0,\overline{h})$. Then, following the proof of Biau, Cadre and
Pellettier [(\citeyear{cadre2007}),
Proposition A.2] if $\varepsilon< \varepsilon_1$, for any $h \in
(0,\overline{h})$,
%
%
\begin{equation}\label{eq:dist.h}
\sup_{x \in\partial L_{h}(\lambda)} \operatorname{dist} \bigl(x,
L_{h}(\lambda+ \varepsilon) \bigr) \leq C_g^{-1} \varepsilon,
\end{equation}
where $C_g$ in the same constant appearing in (\ref{eq:grad})
[see equation (A.1) in \citet{cadre2007}]. In fact, since $C_g$
does not depend on $h$, (\ref{eq:dist.h}) holds uniformly over $h \in
(0, \overline{h})$. Set $\varepsilon(\rho,\tau)= \sup\{ \varepsilon\in
(0,\varepsilon_1) \dvtx C \varepsilon< \tau\rho\}$. Then as $L(\lambda)
\subseteq L_h(\lambda)$ by (C2)(b), (\ref{eq:final}) is verified
for each $\varepsilon< \varepsilon(\rho,\tau)$.
\end{pf}

\begin{appendix}\label{app}
\section*{Appendix: The geometric density}

In this section, we describe in detail our assumptions on the unknown
distribution
$P$. For the sake of completeness, we provide the basic definitions of
Hausdorff measure, Hausdorff dimension, and rectifiability. We refer
the reader to \citet{evansgariepy92}, \citet{mattila99},
\citet{ambrosiofuscopallara00} and \citet{federer69}
for all the relevant geometric and measure theoretic background.

Let $k \in[0,\infty)$. The $k$-dimensional Hausdorff measure of a set
$E$ in $\mathbb{R}^d$ is defined as $\mathcal{H}^k(E) \equiv\lim
_{\delta\downarrow0} \mathcal{H}^k_\delta(E)$, where, for $\delta
\in(0,\infty]$,
\[
\mathcal{H}^k_\delta(E) = \frac{v_k}{2^k} \inf\biggl\{ \sum_{i \in
I} (\operatorname{diam}(E_i))^k \dvtx\operatorname{diam}(E_i) < \delta
\biggr\},
\]
where the infimum is over all the countable covers $\{ E_i \}_{i \in
I}$ of $E$, with the convention $\operatorname{diam}(\varnothing) = 0$.
The Hausdorff dimension of a set $E \subset\mathbb{R}^d$ is
\[
\inf\{ k \geq0 \dvtx\mathcal{H}^k(E) = 0 \}.
\]
Note that $\mathcal{H}^0$ is the counting measure, while $\mathcal
{H}^d$ coincides with the (outer) Lebesgue measure. If $k< d$, we will
refer to any $\mathcal{H}^k$-measureable set as a set of
lower-dimension. When $1 \leq k < d$ is an integer, $\mathcal{H}^k(E)$
coincides with the $k$-dimensional area of $E$, if $E$ is contained in
a $\mathcal{C}^1$ $k$-dimensional manifold embedded in $\mathbb{R}^d$.

The set $E$ is said to be $\mathcal{H}^k$-rectifiable if $k$ is an
integer, $\mathcal{H}^k (E) < \infty$ and there exist countably many
Lipschitz functions $f_i \dvtx\mathbb{R}^k \mapsto\mathbb{R}^d$
such that
\[
\mathcal{H}^k \Biggl( E - \bigcup_{i=0}^\infty f_i(\mathbb{R}^k)
\Biggr).
\]
A Radon measure $\nu$ in $\mathbb{R}^d$ is said to be $k$-rectifiable
if there exists a $\mathcal{H}^k$-rectifiable set $S$ and a Borel
function $f \dvtx S \mapsto\mathbb{R}^d$ such that
\[
\nu(A) = \int_{A \cap S} f(x) \,d \mathcal{H}^k(x)
\]
for each measurable set $A \subseteq\mathbb{R}^d$.

Throughout this article, we assume that $P$ is a finite mixture of
probability measures supported on disjoint compact sets of possibly
different integral Hausdorff dimensions.
Formally, for each Borel set $A \subseteq\mathbb{R}^d$
and for some integer $m$,
\[
P(A) = \sum_{i=1}^m \pi_i P_i(A),
\]
where $\pi$ is a point in the interior of the $(m-1)$-dimensional
standard simplex and each $P_i$ is a $d_i$-rectifiable Radon measure
with compact and connected support $S_i$, where $d_i \in
\{0,1,\ldots,d\}$ and $S_i \cap S_j = \varnothing$,
for each $i \neq j$. Notice that we also have $\max_i \mathcal
{H}^{d_i}(S_i) < \infty$. 
By Theorem 3.2.18 in \citet{federer69}, each of the
lower-dimensional rectifiable sets comprising the support of $P$, can be
represented as the union of $\mathcal{C}^1$ embedded submanifolds,
almost everywhere $P$. Thus, we are essentially allowing $P$ to be a
mixture of distributions supported on disjoint submanifolds of
different dimensions and finite sets.

Our assumptions imply that, for every mixture component $P_i$, there
exists a $\mathcal{H}^{d_i}$-measurable real valued function $p_i$
such that
such that
%
%
\begin{equation}\label{eq:pi}
p_i(x) = \cases{
\displaystyle\lim_{h \to0} \frac{P_i ( B(x,h) )}{v_{d_i} h^{d_i}} >
0, &\quad if $x \in S_i$,\vspace*{2pt}\cr
0, &\quad if $x \notin S_i$,}
\end{equation}
where $v_{d_i}$ is the volume of the unit Euclidean ball in $\mathbb
{R}^{d_i}$. See, for instance,
Mattila [(\citeyear{mattila99}), Corollary 17.9] or Ambrosio, Fusco and
Pallara [(\citeyear{ambrosiofuscopallara00}), Theorem 2.83]. %
Indeed, $p_i$ is a density function with respect to $\mathcal
{H}^{d_i}$ since, for any measurable set~$A$,
\[
P_i(A) = \int_{A \cap S_i} p_i(x) \,d \mathcal{H}^{d_i}(x),
\]
where $\mathcal{H}^{d_i}$ denotes the $d_i$-dimensional Hausdorff
measure on $\mathbb{R}^d$. 

We do not assume any knowledge of the probability measures comprising
the mixture $P$, of their number, supports and dimensions, nor of the
vector of mixing probabilities $\pi$.

Recall that the geometric density is the extended real-valued function
defined as
\[
p(x) = \lim_{h \downarrow0}\frac{ P(B(x,h))}{v_d h^d},\qquad x \in
\mathbb{R}^d.
\]
Below, we list the key properties of the geometric density. Notice, in
particular, that $p$ is not a probability density with respect to $\mu
$, since, in general, $0 \leq\int_{\mathbb{R}^d} p(x) \,d\mu(x) \leq1$.
\begin{proposition}
The geometric density satisfies the following properties:
\begin{longlist}
\item$p(x) = \infty$ if and only if $x \in S_i$ with
$\operatorname{dim}(S_i) < d$, almost everywhere $P$,
\item$p(x) = \pi_i p_i(x) < \infty$ if and only if $x \in
S_i$ with $\operatorname{dim}(S_i) = d$, almost everywhere $\mu$,
\item$\mu( \{ x \dvtx p(x) = \infty\} ) = 0$,
\item if $x \notin S$, then $p(x) = 0$,
\item$S = \overline{\{ x \dvtx p(x) > 0\}}$.
\end{longlist}
\end{proposition}
\begin{pf}
If $S_i$ has dimension $d$, then, by the Lebesgue theorem,
\[
p(x) = \lim_{h \downarrow0}\frac{ P(B(x,h))}{v_d h^d} = \pi_i
p_i(x) < \infty,
\]
$\mu$-almost everywhere on $S_i$. Similarly, if $d_i < d$, then, by
(\ref{eq:pi}),
\[
p(x) = \lim_{h \downarrow0}\frac{ P(B(x,h))}{v_d h^d} = \lim_{h
\downarrow0}\frac{v_{d_i} h^{d_i}}{v_{d} h^{d}}\frac{
P_i(B(x,h))}{v_{d_i} h^{d_i}} = \infty,
\]
since\vspace*{1pt} $\frac{v_{d_i} h^{d_i}}{v_{d} h^{d}} \rightarrow\infty$ as $h
\downarrow0$, $\mathcal{H}^{d_i}$-almost everywhere on $S_i$.
Thus, part (i) and (ii) follow. Part (iii) is a
direct consequence of (i) and (ii), while parts (iv)
and (v) stem directly from the definition of support.
\end{pf}
\end{appendix}

As a final remark, even though the geometric density $p$ is very
different from the mixture densities $p_i$, for our clustering
purposes, we need only to concern ourselves with estimating the level
sets of $p$.

\section*{Acknowledgments}
The authors thank Giovanni Leoni and Aarti Singh for helpful
discussions and two anonymous referees for many constructive comments
and suggestions that greatly improved the exposition.

\printaddresses

\end{document}